\documentclass[11pt]{amsart}

\usepackage{amsmath}
\usepackage{amsxtra}
\usepackage{amscd}
\usepackage{amsthm}
\usepackage{amsfonts}
\usepackage{amssymb}
\usepackage{mathtools}
\usepackage{eucal}
\usepackage[all]{xy}
\usepackage{graphicx}
\usepackage{tikz-cd}
\usepackage{mathrsfs}
\usepackage{euler}
\usepackage{hyperref}
\usepackage{color}
\usepackage{longtable}
\usepackage{float}
\usepackage{caption}
\usepackage{enumerate}
\usepackage{faktor}
\usepackage{bbm}
\usepackage{lineno}
\usepackage{stmaryrd}
\usepackage{rotating}
\usepackage[T1]{fontenc} 
\usepackage[utf8]{inputenc} 

\usepackage[colorinlistoftodos, textsize=tiny]{todonotes}
\def\listtodoname{List of Todos}
\def\listoftodos{\@starttoc{tdo}\listtodoname}

\usepackage{tikz}
\usepackage{tikz-cd}
\usetikzlibrary{matrix,arrows,decorations.pathmorphing}

\usepackage[bottom=1.25in, top=1.5in, right=1.5in, left=1.5in]{geometry}

\theoremstyle{plain}
  \newtheorem{theorem}{Theorem}
	\newtheorem{thm}[theorem]{Theorem}
  \newtheorem{proposition}[theorem]{Proposition}
	\newtheorem{prop}[theorem]{Proposition}
  
  \newtheorem{lem}[theorem]{Lemma}
  \newtheorem{cor}[theorem]{Corollary}
  \newtheorem{corollary}[theorem]{Corollary}

\theoremstyle{definition}
	
	\newtheorem{defn}[theorem]{Definition}
  \newtheorem{example}[theorem]{Example}
	
	\newtheorem{conjecture}[theorem]{Conjecture}
	
\theoremstyle{remark}
	\newtheorem{rem}[theorem]{Remark}
	\newtheorem{remark}[theorem]{Remark}

    {\cleardoublepage\thispagestyle{empty}\null\vfill\begin{center}%
    \bfseries Acknowledgements\end{center}}%

\DeclareMathAlphabet{\mathcal}{OMS}{cmsy}{m}{n}

\newcommand{\one}{\mathbbm{1}}

\newcommand{\C}{\mathbb{C}}

\newcommand{\F}{\mathbb{F}}
\newcommand{\G}{\mathbb{G}}

\newcommand{\Q}{\mathbb{Q}}
\newcommand{\R}{\mathbb{R}}

\newcommand{\Z}{\mathbb{Z}}

\newcommand{\cald}{\mathcal{D}}

\newcommand{\calh}{\mathcal{H}}
\newcommand{\calk}{\mathcal{K}}
\newcommand{\calo}{\mathcal{O}}
\newcommand{\calp}{\mathcal{P}}

\newcommand{\cals}{\mathcal{S}}
\newcommand{\calw}{\mathcal{W}}
\newcommand{\calu}{\mathcal{U}}

\newcommand{\mfs}{\mathfrak{S}}
\newcommand{\rmn}{\mathrm{N}}
\newcommand{\rmo}{\mathrm{O}}

\newcommand{\Ad}{{\mathrm{Ad}}}

\newcommand{\Aut}{{\mathrm{Aut}}}

\newcommand{\Id}{\mathrm{Id}}

\renewcommand{\Im}{{\mathrm{Im}}}

\newcommand{\Ker}{{\mathrm{Ker}}}

\newcommand{\Gal}{{\mathrm{Gal}}}

\newcommand{\Hom}{{\mathrm{Hom}}}

\newcommand{\Ind}{{\mathrm{Ind}}}

\newcommand{\Res}{{\mathrm{Res}}}

\newcommand{\Reg}{\mathrm{Reg}}

\newcommand{\gl}{\mathrm{GL}}
\newcommand{\GL}{{\mathrm{GL}}}
\newcommand{\go}{\mathrm{GO}}
\newcommand{\GO}{{\mathrm{GO}}}

\newcommand{\PSL}{{\mathrm{PSL}}}

\newcommand{\SL}{{\mathrm{SL}}}

\newcommand{\End}{{\mathrm{End}}}

\newcommand{\ord}{{\mathrm{ord}}}
\newcommand{\Tr}{{\mathrm{Tr}}}
\newcommand{\vol}{{\mathrm{vol}}}

\newcommand{\hv}{\mathrm{HV}}

\newcommand{\rs}{\mathrm{RS}}

\newcommand{\new}{{\mathrm{new}}}
\newcommand{\opt}{{\mathrm{opt}}}
\newcommand{\stark}{{\mathrm{Stark}}}

\newcommand{\bs}{\backslash}
\newcommand{\wt}{\widetilde}

\newcommand{\longhookrightarrow}{\, \xhookrightarrow{\quad} \,}
\newcommand{\lra}{{\, \longrightarrow \,}}
\newcommand{\iso}{\, \xrightarrow{\widesim{}} \,}

\newcommand{\paren}[1]{\mathopen{}\left(#1\right)\mathclose{}}
\newcommand{\set}[1]{\mathopen{}\left\{#1\right\}\mathclose{}}
\newcommand{\sbrac}[1]{\mathopen{}\left[#1\right]\mathclose{}}
\newcommand{\abrac}[1]{\mathopen{}\left\langle#1\right\rangle\mathclose{}}
\newcommand{\verts}[1]{\mathopen{}\left\lvert#1\right\rvert\mathclose{}}
\newcommand{\norm}[1]{\mathopen{}\left\lvert\left\lvert#1\right\rvert\right\rvert\mathclose{}}
\newcommand\restr[2]{{
  \left.\kern-\nulldelimiterspace 
  #1 
  \right|_{#2} 
  }}

\newcommand{\Mid}{\,\middle|\,}

\newcommand{\pair}[1]{\abrac{#1}}
\newcommand{\wh}{\widehat}

\newcommand{\widesim}[2][2]{
  \mathrel{\overset{#2}{\scalebox{#1}[1]{$\sim$}}}
}

\renewcommand{\setminus}{-}

\makeatletter
\AtBeginDocument
 {
    \def\@thm#1#2#3{%
      \ifhmode
        \unskip\unskip\par
      \fi
      \normalfont
      \trivlist
      \let\thmheadnl\relax
      \let\thm@swap\@gobble
      \let\thm@indent\indent 
      \thm@headfont{\scshape}
      \thm@notefont{\fontseries\mddefault\upshape}%
      \thm@headpunct{.}
      \thm@headsep 5\p@ plus\p@ minus\p@\relax
      \thm@space@setup
      #1
      \@topsep \thm@preskip               
      \@topsepadd \thm@postskip           
      \def\dth@counter{#2}%
      \ifx\@empty\dth@counter
        \def\@tempa{%
          \@oparg{\@begintheorem{#3}{}}[]%
        }%
      \else
        \H@refstepcounter{#2}%
        \hyper@makecurrent{#2}%
        \let\Hy@dth@currentHref\@currentHref
        \AddToHookNext{para/begin}{\MakeLinkTarget*{\Hy@dth@currentHref}}%
        \def\@tempa{%
          \@oparg{\@begintheorem{#3}{\csname the#2\endcsname}}[]%
        }%
      \fi
      \@tempa
    }%
  \dth@everypar={%
    \@minipagefalse \global\@newlistfalse
    \@noparitemfalse
    \if@inlabel
      \global\@inlabelfalse
      \begingroup \setbox\z@\lastbox
       \ifvoid\z@ \kern-\itemindent \fi
      \endgroup
      \unhbox\@labels
    \fi
    \if@nobreak \@nobreakfalse \clubpenalty\@M
    \else \clubpenalty\@clubpenalty \everypar{}%
    \fi
  }%
}

\title[The Harris--Venkatesh conjecture for derived Hecke operators III]
{The Harris--Venkatesh conjecture for derived Hecke operators III:
	local constants}

\author{Robin Zhang}
\address{Department of Mathematics, Columbia University, USA \newline
	\indent Department of Mathematics, Massachusetts Institute of Technology, USA \newline
  \indent Present address: Institut des Hautes \'{E}tudes Scientifiques, France}
\email{rzhang@ihes.fr}

\date{August 25, 2026}

\begin{document}

\begin{abstract}
	\normalsize
	The first two papers in this series prove
	the Harris--Venkatesh conjecture on derived Hecke operators
	and its refinement with Stark units
	for imaginary dihedral modular forms of weight $1$.
	This paper explicitly describes the
	constants appearing in the
	Harris--Venkatesh conjecture
	for dihedral modular forms
	by evaluating $\mathrm{GL}(2) \times \mathrm{GL}(2)$
	Rankin--Selberg periods and zeta integrals
	on newforms and optimal forms. One consequence
	is a formula for the ratio between
	Petersson norms and adjoint $L$-values.
	These calculations also extend to
	exotic modular forms of odd level
	and to exotic modular forms with $2$-ordinary
	Deligne--Serre representation.
\end{abstract}

\maketitle

\setcounter{tocdepth}{1}
\tableofcontents


\section{Introduction}
Let $f$ be any Hecke newform of weight $1$
for the congruence subgroup $\Gamma_1(N)$,
and let $\rho_f$ be its associated Artin representation.
Choose a finite Galois extension
$E/\Q$ through which $\Ad(\rho_f)$ factors, and set
$G := \Gal(E/\Q)$ and
$R_f := \Z[\chi_{\Ad(\rho_f)}]$.
Let $V_{\Ad(\rho_f)}$ be the space of complex
$2 \times 2$ matrices with trace $0$, with the action
of $G$ by conjugation through $\rho_f$. Let
$M_{\Ad} \subset V_{\Ad(\rho_f)}$ be the $G$-stable
$R_f$-lattice generated by
\[
	x_\sigma
		:= 2\rho_f(\sigma)
			- \Tr\paren{\rho_f(\sigma)}
				\Id_{2 \times 2},
	\qquad \sigma \in G.
\]
As in \cite[Introduction]{zhang-hv},
define the dual unit group
\[
	\calu\paren{\Ad(\rho_f)}
		:= \Hom_{R_f[G]}
			\paren{
				M_{\Ad},
				\calo_E^\times \otimes_\Z R_f}.
\]
Letting $\Z_{(p - 1)}$ be the localization of $\Z$
at the integers coprime to $p - 1$, there are regulator
maps
\begin{align*}
	\Reg_\R:
	\calu\paren{\Ad(\rho_f)} \otimes_\Z \Q
		&\lra \R \otimes_\Z R_f, \\
	\Reg_{\F_p^\times}:
	\calu\paren{\Ad(\rho_f)}
		\otimes_\Z \Z_{(p - 1)}
		&\lra \F_p^\times
			\otimes_\Z R_f\sbrac{\frac{1}{6N}}.
\end{align*}
$\Reg_{\R}$ is the Stark regulator map,
whose relation to the leading coefficient $L'(\Ad(\rho_f), 0)$
of the Artin $L$-function is the subject of the Stark conjecture
\cite{stark-1971,stark-1975,stark-1976,stark-1980}.
The compatibility of the regulator map
$\Reg_{\F_p^\times}$ with the action of derived
Hecke operators on $f$ is the subject of the
Harris--Venkatesh conjecture.

Let $f^* = \sum_n \overline{a_n}q^n$ be the dual
newform of $f = \sum_n a_nq^n$, so
$f^*(z) = \overline{f(-\overline{z})}$.
For every prime $p \nmid 6N$, let $\mfs_p$ denote
the Shimura-class functional of
\cite[Section~3.1]{hv}. Following
\cite[Section~1]{zhang-hvs}, define the
Harris--Venkatesh norm by
\[
	\norm{f}^2_{\F_p^\times}
		:=
		\mfs_p\paren{
			\Tr_p^{Np}\paren{f(z)f^*(pz)}}
		\in
		\F_p^\times
		\otimes
		R_f\sbrac{\frac{1}{6N}}.
\]
In this notation, the conjecture of
Harris--Venkatesh
\cite[Conjecture~3.1]{hv}
takes the following form.
\begin{conjecture}[Harris--Venkatesh]
	\label{conj:hv}
	Let $f$ be a Hecke new cusp form of weight $1$
	and level $N$.
	There is an element
	$u \in \calu\paren{\Ad(\rho_f)}$
	and a positive integer $m$ such that
	\[
			m \cdot \norm{f}^2_{\F_p^\times}
				= \Reg_{\F_p^\times}\paren{u}
		\]
	for every prime $p \nmid 6N$.
\end{conjecture}

The first paper in this series
\cite{zhang-hv} proves
Conjecture~\ref{conj:hv} for all imaginary
dihedral forms, extending earlier results of
Darmon--Harris--Rotger--Venkatesh
\cite{dhrv} and Lecouturier
\cite{lecouturier-hv}.
\begin{theorem}[{\cite[Theorem~4]{zhang-hv}}]
	\label{thm:zhang-hv}
	Let $f$ be an imaginary dihedral modular form
	of weight $1$.
	Then Conjecture~\ref{conj:hv} for $f$ is true.
\end{theorem}

The second paper in this series formulates a unified
conjecture combining the Harris--Venkatesh and
Stark conjectures and requiring the same unit
to satisfy both regulator identities
\cite[Conjecture~6]{zhang-hvs}.
It refines both the
Harris--Venkatesh conjecture
and the $F = \Q$ case of
Horawa's conjecture
\cite[Conjecture~1.1]{horawa},
and proves the combined statement in the
imaginary dihedral case.
\begin{theorem}[{\cite[Theorem~1]{zhang-hvs}}]
	\label{thm:zhang-hvs}
	Let $f$ be an imaginary dihedral newform
	of weight $1$ and level $N$.
	There is a unique
	(up to multiplication by roots of unity)
	$u_f \in
	\calu\paren{\Ad(\rho_f)} \otimes \Q$
	such that
	\begin{enumerate}[(a)]
		\item
			\[
				\norm{f}^2_{\R}
					= \Reg_\R\paren{u_f},
			\]
			compatible with conjugations of $f$
			under $\Aut(\C)$;
		\item for every sufficiently large prime
			$p \nmid 6N$,
			\[
				\norm{f}^2_{\F_p^\times}
					= \Reg_{\F_p^\times}\paren{u_f}.
			\]
	\end{enumerate}
\end{theorem}
One of the main goals of this third paper is to
calculate the local Rankin--Selberg scalar comparing
the newform and optimal form. In the imaginary
dihedral case, this calculation makes the data in
Conjecture~\ref{conj:hv} explicit: one may take
$m = m_f$ and $u = m_f u_f$ for an explicitly
computable positive integer $m_f$.
The other main goal is to extend the theory of
optimal forms, introduced for imaginary dihedral
forms in \cite[Definition~6.1]{zhang-hv},
to as many exotic newforms as possible.

\subsection*{Rankin--Selberg comparison for
locally dihedral forms}
Constructing this multiplier requires a detailed
local theory of the Rankin--Selberg period
defined globally in \cite[Section~4]{zhang-hvs}.
We prove our main local result more generally than
the Harris--Venkatesh and Stark conjectures are known,
since the local calculation uniformly treats dihedral newforms
and exotic newforms that are locally dihedral.

Outside of the numerical evidence given by
Marcil \cite{marcil}, all of the current progress on the
Harris--Venkatesh conjecture concerns dihedral forms.
Our local calculations apply to dihedral forms
and a more general class of modular forms of weight $1$:
we call a newform $f$ \emph{locally dihedral} if,
at every prime $p \mid N$,
its local representation $\rho_{f, p}$
is either induced from a character of
a quadratic field extension $K_p / \Q_p$
or is a sum of two characters.
In particular, $f$ is locally dihedral
if $N$ is odd or if $\rho_f$ is $2$-ordinary.
For a locally dihedral form, the characters
$\chi_N := \prod_{p \mid N} \chi_p$ and
$\xi_N := \prod_{p \mid N} \xi_p$
and the conductor exponents $o(\xi_p)$
described in Section~\ref{sec:locally-dihedral}
will play the same role as
$\chi$, $\xi$, and $o(\xi)$ do for dihedral forms.

For the automorphic representation $\pi_f$ generated
by $f$ and its contragredient
$\wt{\pi}_f = \pi_{f^*}$,
\cite[Section~4]{zhang-hvs} defines the global
Rankin--Selberg period by
\[
	\calp_\rs\paren{f_1 \otimes f_2}
		:= \frac{1}{\sbrac{\PSL_2(\Z) : \Gamma}}
			\int_{X(\Gamma)}
				f_1(z)f_2(-\overline{z})y\,d\mu.
\]
The Harris--Venkatesh period $\calp_\hv$ is defined
through the Shimura class in
\cite[Section~2]{zhang-hv}.
For $f^\new := f \otimes f^*$, one has
\begin{align*}
	\calp_\rs\paren{f^\new}
		&= \sbrac{\PSL_2(\Z) : \Gamma_0(N)}^{-1}
			\norm{f}_\R^2, \\
	\calp_\hv\paren{f^\new}
		&= \sbrac{\PSL_2(\Z) : \Gamma_0(N)}^{-1}
			\norm{f}_{\F_p^\times}^2.
\end{align*}

For imaginary dihedral $f$, let $\varphi^\opt$ be
the automorphic avatar of the optimal form.
An element $u_{\varphi^\opt}$ satisfying
the relations
\begin{align*}
	\calp_\rs\paren{\varphi^\opt}
		&= \Reg_\R\paren{u_{\varphi^\opt}} \\
	\calp_\hv\paren{\varphi^\opt}
		&= \Reg_{\F_p^\times}
			\paren{u_{\varphi^\opt}}
\end{align*}
is constructed in \cite[Lemma~5.2]{zhang-hvs}.
The ratio of the Rankin--Selberg periods on the
newform and optimal form is the characteristic-$0$
scalar entering the integral comparison of
\cite[Section~4]{zhang-hv} giving
a cross-multiplied relation between
the corresponding Harris--Venkatesh periods;
we calculate these local ratios explicitly.
\begin{theorem}
	\label{thm:ratio}
	If $f$ is a locally dihedral newform of weight $1$ and level $N$
	with associated character $\chi_N$, then
	$\frac{\calp_\rs\paren{f^\new}}{\calp_\rs\paren{f^\opt}}
	\in \Q(\xi_N + \xi_N^{-1})^\times$
	and
	\[
		\frac{\calp_\rs\paren{f^\new}}{\calp_\rs\paren{f^\opt}}
			= \prod_{p \mid N} \alpha_{\chi_p},
	\]
	where $\alpha_{\chi_p}$ depends only on $\chi_p$.
	Moreover, for odd primes $p$ that are not simultaneously ramified
	in both $K_p$ and $\chi_p$:
	\begin{enumerate}
		\item If $p$ is ramified in $K_p$, then 
			\[
				\alpha_{\chi_p} = 4.
			\]
		\item If $p$ is inert in $K_p$,  then 
			\[
				\alpha_{\chi_p} =
					\begin{cases}
						\frac{(p - 1)^2}{p^2}
							&\text{if $\xi_p$ is unramified,} \\
						\frac {\xi_p(-1)(p + 1)^2}{p^2}
							&\text{if $\xi_p$ is ramified and $\xi_p^2$ is unramified,}\\
						\frac {\xi_p(-1)(p + 1)}{p^2}
							&\text{if $\xi_p^2$ is ramified.}
					\end{cases}
			\]
		\item If $p$ is split in $K_p$ (so $K_p = \Q_p \oplus \Q_p$ is split
			with uniformizers $\varpi_1, \varpi_2$
			and $\chi_p = (\chi_1, \chi_2)$), then
			\[
				\alpha_{\chi_p} =
					\begin{cases}
						\frac{(p - 1)(p - \xi_p(\varpi_1))
						(p - \xi_p(\varpi_2))}{(p + 1) p^2}
							&\text{if $\chi_p$ is ramified and $\xi_p$ is unramified,} \\
						\frac{\chi_p(-1) (p - 1)^2 p^{2o(\xi_p) - 1}}
						{p^{2o(\xi_p) + 1} - p^{2o(\xi_p)} + 2}
							&\text{if $\xi_p$ is ramified, $\xi_p^2$ is unramified,} \\
							&\text{and exactly one of the $\chi_i$ is ramified,} \\
						\frac{\chi_p(-1) (p - 1)^3 p^{2o(\xi_p) - 2}}
						{p^{2o(\xi_p) + 1} - p^{2o(\xi_p)} + 2}
							&\text{if $\xi_p$ is ramified, $\xi_p^2$ is unramified,} \\
							&\text{and both $\chi_i$ are ramified,} \\
						\frac{\chi_p(-1) (p - 1)^2 p^{2o(\xi_p) - 1}}
						{(p + 1) \paren{p^{2o(\xi_p) + 1} - p^{2o(\xi_p)} + 2}}
							&\text{if $\xi_p^2$ is ramified} \\
							&\text{and exactly one of the $\chi_i$ is ramified,} \\
						\frac{\chi_p(-1) (p - 1)^3 p^{2o(\xi_p) - 2}}
						{(p + 1)(p^{2o(\xi_p) + 1} - p^{2o(\xi_p)} + 2)}
							&\text{if $\xi_p^2$, $\chi_1$, and $\chi_2$ are all ramified.}
					\end{cases}
			\]
	\end{enumerate}
\end{theorem}


\subsection*{Petersson norm constants}
Stark \cite[Section 6]{stark-1975}
observed that the derivative $L'(\Ad(\rho_f), 0)$ of the Artin $L$-function
at $s = 0$ (alternatively, the value $L(\Ad(\rho_f), 1)$),
is related to the Petersson norm $\norm{f}_\R^2$
by a constant $c_{f, \rs}$ involving a finite product of Euler factors
depending only on $f$.
We prove a precise formula for $c_{f, \rs}$.
\begin{proposition}
	\label{prop:crs}
	There is a $c_{f, \rs} \in \Q(\chi_{\rho_f})$ such that
	\[
		\norm{f}_\R^2 = c_{f, \rs} \cdot L'\big(\Ad(\rho_f), 0 \big).
	\]
	Moreover in terms of the Euler factors $\zeta_p(s)$,
	$L_p(f \times f^*, s)$, and $L_p(\Ad(\rho_f), s)$
	of the Riemann zeta function, Rankin--Selberg $L$-function,
	and Artin $L$-function,
	\[
		c_{f, \rs} = \restr{\paren{\prod_{p \mid N}
			\frac{L_p\paren{f \times f^*, s}}
			{\zeta_p(s) L_p\big(\Ad(\rho_f), s\big)}}}
			{s = 0} \, .
	\]
\end{proposition}

Combining Theorem~\ref{thm:ratio} with
\cite[Theorem~1]{zhang-hvs} gives the following
description of $c_{f, \rs}$ in the imaginary
dihedral case.
Recall that for each imaginary dihedral cusp form $f$
there is an imaginary quadratic number field $K$
and a character $\chi$ of $G_K$ whose induction to
$G_{\Q}$ is isomorphic to the Deligne--Serre
representation $\rho_f$.
Let $c_K$ denote the nontrivial element of
$\Gal(K/\Q)$. The associated antinorm
$\xi := \chi^{1 - c_K}$ has conductor
$c = c(\xi)$ and ring class field $H_c$.
Let $H_1$ be the Hilbert class field of $K$, and let
$w_K$ be the number of roots of unity in $K$.
\begin{corollary}
	\label{cor:crs}
	Let $f$ be an imaginary dihedral cusp form of weight $1$
	and level $\Gamma_1(N)$,
	$K$ be the corresponding imaginary quadratic number field,
	and $\chi$ be the corresponding character of $G_K$.
	Then
	\[
		c_{f, \rs}
			= - \frac{\sbrac{H_c : H_1}w_K}{2}
				\prod_{p \mid N} \alpha_{\chi_p},
	\]
	where $\alpha_{\chi_p}$ is as described in
	Theorem \ref{thm:ratio}.
\end{corollary}

\subsection*{Harris--Venkatesh constants}
Finally, let us return to the integer $m$ in
Conjecture~\ref{conj:hv}.
Let $D_K$ be the discriminant of $K/\Q$,
let $\delta_{K/\Q}$ be its different,
let $c(\chi)$ be the conductor of $\chi$, and let
$\verts{\Im(\chi^2)}$ denote the cardinality of the
image of $\chi^2$.
Experimental data on $m$ have been compiled in
dihedral cases by \cite{hv, marcil},
and theoretical upper bounds are known
due to \cite[Remark~1.3]{dhrv}
in the following cases:
\begin{itemize}
	\item one may choose $m$ dividing $24$ if $f$ is
		real dihedral, $D_K$ is odd, and
		$c(\chi) \mid \delta_{K/\Q}$;
	\item one may choose $m$ dividing $6q$ if $f$ is
		imaginary dihedral, $D_K$ is an odd prime,
		$\chi$ is unramified, and
		$\verts{\Im(\chi^2)}$ is a power of a prime $q$;
	\item one may choose $m$ dividing $6$ if $f$ is
		imaginary dihedral, $D_K$ is an odd prime,
		$\chi$ is unramified, and
		$\verts{\Im(\chi^2)}$ is not a prime power.
\end{itemize}

For imaginary dihedral $f$,
Corollary~\ref{cor:crs} and the uniqueness in
Theorem~\ref{thm:zhang-hvs} identify the canonical
unit as $u_f = c_{f, \rs} u_\stark$, where
$c_{f, \rs}$ is expressed in terms of the local
constants $\alpha_{\chi_q}$.
The integral comparison of
\cite[Section~4]{zhang-hv} makes the corresponding
multiplier explicit.

\begin{corollary}
	\label{cor:m}
	Let $f$ be an imaginary dihedral newform of weight $1$
	and level $\Gamma_1(N)$.
	Let $m_f$ be the explicit positive integer defined in
	Equation~\ref{eq:explicit-multiplier}.
	Then in Conjecture~\ref{conj:hv}, one may take
	$m = m_f$ and
	$u = m_f u_f = m_f c_{f, \rs} u_\stark$.
\end{corollary}

\begin{remark}
	For any locally dihedral form $g$ whose
	corresponding optimal form $g^\opt$ is known
	to satisfy compatible regulator identities,
	the same construction would give explicit
	descriptions of $c_{g, \rs}$ and $m_g$.
	We expect this to apply to real dihedral forms
	in the forthcoming work \cite{zhang-real}.
\end{remark}

\subsection*{Outline}
Section~\ref{sec:local-background} develops the
theory of Rankin--Selberg periods and extends the
definition of optimal forms to locally dihedral forms.
Section~\ref{sec:non-vanishing} proves their
non-vanishing and rationality properties at new
vectors and optimal vectors.
Section~\ref{sec:inner-product-opt} explicitly
calculates Rankin--Selberg zeta integrals evaluated
at optimal vectors.
Section~\ref{sec:inner-product-new} explicitly
calculates Rankin--Selberg zeta integrals evaluated
at new vectors.
Section~\ref{sec:calculation-local} combines these
calculations to prove Theorem~\ref{thm:ratio}.
Finally, Section~\ref{sec:petersson} proves
Proposition~\ref{prop:crs} and deduces
Corollaries~\ref{cor:crs} and~\ref{cor:m}.

\subsection{Acknowledgements}

I am deeply thankful to
Michael Harris, who supervised the thesis of which
this article is roughly the second half,
for suggesting this area of research,
for sharing his deep and broad insight,
and for his continued guidance throughout my doctoral studies.
I am grateful to Masao Oi and the
Department of Mathematics at Kyoto University for their
incredible hospitality and for creating a great environment
for the final stages of preparing this manuscript.
Finally, I would like to thank Wee-Teck Gan
for helpful discussions
and the anonymous referee for making careful suggestions to improve this exposition.

This material is based on work
supported by the National Science Foundation
under Grants DGE-1644869, DGE-2036197, and DMS-2303280.


\numberwithin{equation}{section}
\numberwithin{theorem}{section}


\section{Local Rankin--Selberg periods and Whittaker functions}
\label{sec:local-background}

In Sections \ref{sec:local-essential}
and \ref{sec:local-theta-liftings},
we recall non-archimedean local background
(with uniform notation for our setting) from:
\begin{itemize}
	\item $\GL_2$ theory following
		Jacquet--Langlands~\cite{jacquet-langlands}, Casselman \cite{casselman},
		Bump \cite{bump}, and Zhang \cite{zhang-asian},
	\item Rankin--Selberg $\GL_2 \times \GL_2$ theory following Jacquet~\cite{jacquet}, and
	\item $\GL_2 \times_{\G_m} \GO(V)$ theta liftings following
		Waldspurger \cite{waldspurger} and Yuan--Zhang--Zhang \cite{yzz}).
\end{itemize}
We use these tools to extend the local theory of Rankin--Selberg periods
for newforms and optimal forms
in Sections \ref{sec:quadratic-cases} and \ref{sec:locally-dihedral}.

Throughout this article, let $p$ be any prime,
let $F$ be a $p$-adic field
with ring of integers $\calo_F$ and residue field $\F_q$,
and fix a non-trivial additive character
$\psi: F \lra \C^\times$.
Let $E/F$ be a semisimple $F$-algebra of
dimension $2$ with trace map $\Tr: E \lra F$ and norm map
$\rmn: E^\times \lra F^\times$, and let
$\chi$ be a character of $E^\times$.
Using the Weil representation, there is an irreducible representation
$\pi(\chi)$ of $\GL_2(F)$ with central character
$\omega = \eta \cdot \restr{\chi}{F^\times}$, where $\eta$ is a
quadratic character on $F^\times$ with kernel $\rmn(E^\times)$.

\subsection{Local Rankin--Selberg periods}
\label{sec:local-essential}

\subsubsection{Local Whittaker models and Kirillov models for $\texorpdfstring{\gl_2}{GL(2)}$}
\label{subsec:whittaker-kirillov}

The induced representation
$\calw(\psi) := \Ind_{N(F)}^{\GL_2(F)}(\psi)$ is 
the space of functions $W$ on $\GL_2(F)$ such that
\[
	W\paren{\begin{pmatrix}1 & x \\ & 1 \end{pmatrix} g} = \psi(x) W(g),
\]
with action by $\GL_2(F)$ via translation.

Let $\pi$ be an irreducible infinite-dimensional representation of
$\GL_2(F)$. We can embed $\pi$ into $\calw(\psi)$:
\[
	\pi \longhookrightarrow \calw(\psi).
\]
This embedding is unique up to scaling
(cf. \cite[Theorem~4.1.2]{bump}).
Therefore, we have a well-defined subspace
$\calw(\pi, \psi) \subset \calw(\psi)$ called the Whittaker model of
$\pi$.
If we change $\psi$ to
$\psi_a(x) := \psi(ax)$, then we have the isomorphism
\[
	\begin{tikzcd}[row sep = tiny]
		\calw(\psi) \arrow[r, "\widesim{}"]
			& \calw\paren{\psi_a} \\
		W(g) \arrow[r, mapsto]
			&
		W\paren{
			\begin{pmatrix} a & \\ & 1 \end{pmatrix}g}
	\end{tikzcd}
\]
Thus, without loss of generality, we can assume that
$\psi$ has order $0$ in the sense that $\calo_F$ is
the maximal fractional ideal of $F$ over which
$\psi$ is trivial.

Define
\[
	\begin{tikzcd}[row sep = tiny]
		F^\times \arrow[r, "a"]
			& \GL_2(F) \\
		x \arrow[r, mapsto]
			& \begin{pmatrix} x & 0 \\ 0 & 1 \end{pmatrix}
	\end{tikzcd}
\]
According to Jacquet--Langlands
\cite[Section~2]{jacquet-langlands}
(cf. \cite[Section~4.4]{bump}),
the following restriction map is injective,
\[
	W(g) \longmapsto \kappa(x) := W\big(a(x)\big).
\]
Let $\calk (\pi, \psi)$ be the image of this map with the induced action
by $\GL_2(F)$. This is called the Kirillov model.
The action of the Borel subgroup of $\GL_2(F)$
on $\calk(\pi, \psi)$ is as follows
(cf. \cite[Section~14]{jacquet}, \cite[Equation~4.25]{bump}),
\[
	\begin{pmatrix} a & b \\ 0 & d \end{pmatrix} \kappa(x)
		= \omega(d) \psi\paren{\frac{bx}{d}} \kappa\paren{\frac{ax}{d}}.
\]
The action of $w := \begin{psmallmatrix} 0 & 1 \\ -1 & 0 \end{psmallmatrix}$
is not easy to write down (cf. \cite[Section~4.7]{bump}). 

If $\pi$ is changed to $\pi \otimes \mu$ for a character $\mu$ of
$F^\times$ via $\det: \GL_2(F) \lra F^\times$, then 
$\calw(\pi \otimes \mu, \psi) = \calw(\pi, \psi) \otimes \mu$.

Let $\varpi$ be a uniformizer of $F$.
For any non-negative integer $i$, define $U_0(\varpi^i)$ and $U_1(\varpi^i)$
as the following subgroups of $\GL_2(\calo_F)$
(cf. \cite[Section~2.3]{zhang-asian}),
\begin{align*}
	U_0(\varpi^i) :&= \set{\begin{pmatrix} a & b \\ c & d \end{pmatrix}
		\Mid c\equiv 0 \pmod {\varpi^i}}, \\
	U_1(\varpi^i) :&= \set{\begin{pmatrix} a & b \\ c & d \end{pmatrix}
		\Mid (c, d) \equiv (0, 1)\pmod {\varpi^i}}.
\end{align*}
For any irreducible representation $\pi$ of
$\GL_2(F)$, define its level $o = o(\pi)$ to be the
least non-negative integer $i$ such that
$\pi^{U_1(\varpi^i)} \neq 0$.
Let $\calw(\pi) = \calw(\pi,\psi)$ be the
Whittaker model of $\pi$.
The space
$\calw(\pi)^{U_1(\varpi^o)}$ is one-dimensional
\cite[Section~2.3]{zhang-asian}.
Let $W^\new$ be its unique element satisfying
$W^\new(e) = 1$.
If $\psi$ is changed to $\psi_a$ for some
$a \in \calo_F^\times$, then $W^\new(g)$ is
changed to
$W^\new\paren{
	\begin{pmatrix} a & \\ & 1 \end{pmatrix}g}$.
There is no simple formula for the
change in $W^\new$ when $\pi$ is replaced by
$\pi\otimes\mu$.

Following Jacquet--Langlands~\cite[Theorem~2.18]{jacquet-langlands}
(cf. \cite[Section~14]{jacquet}, \cite[Proposition~4.7.5]{bump}),
we have the zeta integral
\[
	Z(W, s) := \int_{F^\times}
		W\paren{\begin{pmatrix} a & \\ & 1 \end{pmatrix}}
		\verts{a}^{s - \frac{1}{2}} da,
\]
where $da$ is a Haar measure on $F^\times$ normalized such that
$\vol(\calo_F^\times) = 1$.
These integrals are absolutely convergent when $\Re(s) \gg 0$.
The values of these integrals define a fractional ideal of
$\C[q^{\pm s}]$ with a generator
\[
	L(\pi, s) = Z\paren{W^\new, s}.
\]
Then we can define the normalized zeta integral
\[
	\Psi(W, s) := \frac{Z(W, s)}{L(\pi, s)} \in \C\sbrac{q^{\pm s}}.
\]

For example, if $\pi = \pi(\chi_1, \chi_2)$ is a principal series with
$\chi_1, \chi_2$ unramified and $\alpha_i = \chi_i(\varpi)$
(cf. \cite[Section~3]{jacquet-langlands}),
the above identity is equivalent to
\[
	\frac{1}{\paren{1 - \alpha_1 q^{-s}}\paren{1 - \alpha_2 q^{-s}}}
		= \sum_{n = 0}^\infty W^\new \paren{
			\begin{pmatrix} \varpi^n & 0 \\ 0 & 1 \end{pmatrix}}
			q^{-n \paren{s - \frac{1}{2}}}.
\]
It follows (cf. \cite[Proposition 3.5]{jacquet-langlands})
that $W^\new$ is the unique element in $\calw(\psi)$
which is invariant under $\GL_2(\calo_F)$
with central character $\chi_1 \chi_2$ and takes values
\[
	W^\new \paren{\begin{pmatrix} \varpi^n & 0 \\ 0 & 1 \end{pmatrix}}
		= q^{\frac{-n}{2}}
		\frac{\alpha_1^{n + 1} - \alpha_2^{n + 1}}{\alpha_1 - \alpha_2}.
\]

\subsubsection{Local Rankin--Selberg periods for
	\texorpdfstring{$\GL_2 \times \GL_2$}{GL(2) x GL(2)}}
\label{sec:rs-pairing}

Let $\pi$ be an irreducible infinite-dimensional representation of
$\GL_2(F)$ with contragredient representation $\wt{\pi}$. 
Then we have a canonical $\GL_2(F)$-invariant pairing:
\[
	\calp: \pi \otimes \wt{\pi} \lra \C.
\]
We realize the representations $\pi$ and $\wt{\pi}$
with respective Whittaker models
$\calw(\pi, \psi)$ and $\calw(\wt\pi, \psi)$. 
We want to explicitly construct a $\GL_2(F)$-invariant pairing,
\[
	\calp_{\rs, \psi}: \calw(\pi, \psi) \otimes \calw(\wt{\pi}, \psi)
		\lra \C,
\]
using the Rankin--Selberg method.
This pairing will induce a unique isomorphism with a compatible
linear form:
\[
	\pi \otimes \wt{\pi} \iso
		\calw(\pi, \psi) \otimes \calw(\wt{\pi}, \psi).
\]

First, consider the right-hand side as a space of
functions on
$\GL_2(F) \times \GL_2(F)$. We then restrict this
space to the diagonal:
\[
	\begin{tikzcd}[row sep = tiny]
		\calw(\pi,\psi)
			\otimes\calw(\wt{\pi},\psi)
			\arrow[r, "\Delta^*"]
			&
		\Ind_{N(F)Z(F)}^{\GL_2(F)}(\one) \\
		W_1\paren{g_1}\otimes W_2\paren{g_2}
			\arrow[r, mapsto]
			& W_1(g)W_2(\epsilon'g)
	\end{tikzcd}
\]
where
$\epsilon'
:= \begin{psmallmatrix} 1 & 0 \\ 0 & -1 \end{psmallmatrix}$.
We use the Iwasawa decomposition
\[
	\GL_2(F) = B(F) \GL_2(\calo_F)
\]
to define a function $f(g, s)$ on $N(F)Z(F) \bs \GL_2(F)\times \C$ by
\[
	f\paren{\begin{pmatrix} a & x \\ & b \end{pmatrix} k, s}
		=	\verts{\frac{a}{b}}^s,
\]
for $k \in \GL_2(\calo_F)$.

Following Jacquet~\cite[Section~14]{jacquet},
we consider the zeta integral
\begin{equation}
	\label{eq:zeta-integral-rs}
	Z\paren{W_1, W_2, s} := \int_{N(F)Z(F) \bs \GL_2(F)}
		W_1(g) W_2(\epsilon' g) f(g, s) dg.
\end{equation}
Here $dg$ is the quotient measure on $N(F) Z(F) \bs \GL_2(F)$
constructed from the Haar measures
on $\GL_2(F)$ and its various subgroups.
Write
\[
	g = 
		\begin{pmatrix} 1 & x \\ & 1 \end{pmatrix}
		\begin{pmatrix} z & \\ & z \end{pmatrix}
		\begin{pmatrix} a & \\ & 1 \end{pmatrix}
		k,
\]
with $a, x, z \in F$ and $k \in \GL_2(\calo_F)$.
Then
\[
	dg = dx dz \frac{da}{\verts{a}} dk,
\]
where $dx$ is a measure on $F$ so that $\vol(\calo_F) = 1$,
$dz$ and $da$ are measures on $F^\times$ so that
$\vol(\calo_F^\times) = 1$,
and $dk$ is a measure on $\GL_2(\calo_F)$ with $\vol(\GL_2(\calo_F)) = 1$.

\begin{remark}
\label{rem:epsilon}
	Note that we use
	$\epsilon' := \begin{psmallmatrix} 1 & 0 \\ 0 & -1 \end{psmallmatrix}$
	instead of
	$\epsilon := \begin{psmallmatrix} -1 & 0 \\ 0 & 1 \end{psmallmatrix}$
	in order to notationally avoid the repeated appearance of
	$\omega(-1)$ in our calculations.
	$\epsilon$ is called $\eta$ by Jacquet~\cite[p.11]{jacquet}
\end{remark}

The above zeta integral is absolutely convergent when $\Re(s) \geq 0$,
and has values forming a fractional ideal of $\C[q^{\pm s}]$ with
generator,
\[
	(1 + q^{-s}) L\big(\Ad(\pi), s\big).
\]
We can then define the normalized zeta integral
(cf. \cite[Theorem~14.7]{jacquet}):
\begin{equation}
	\label{eq:normalized-zeta-rs}
	\Psi\paren{W_1, W_2, s}
		= \frac{Z\paren{W_1, W_2, s}}
		{\paren{1 + q^{-s}} L\big(\Ad(\pi), s\big)}
		\in \C\sbrac{q^{\pm s}}.
\end{equation}
We define the invariant form $\calp_{\rs,\psi}$ by
\[
	\calp_{\rs,\psi}\paren{W_1 \otimes W_2}
		:= Z\paren{W_1, W_2, 1}.
\]
Note that up to an $\epsilon$-factor of $\Ad(\pi)$,
we may replace $\Psi(W_1, W_2, 1)$ by $\Psi(W_1, W_2, 0)$,
which is a regularization for the usual divergent integral,
\[
	\int_{N(F)Z(F) \bs \GL_2(F)} W_1(g) W_2(\epsilon' g) dg.
\]

Note that if we change $\psi$ to $\psi_a$, then the above pairings
are compatible with the isomorphism
\[
	\calw(\pi, \psi) \otimes \calw(\wt{\pi}, \psi)
		\iso \calw(\pi, \psi_a) \otimes \calw(\wt{\pi}, \psi_a).
\]
So we may drop the $\psi$ subscript from $\calp_{\rs,\psi}$.
Also, we note that the form $W^\new \otimes \wt{W}^\new$ is invariant
under the canonical isomorphisms with respect to different $\psi$'s.

\begin{example}
	\label{ex:rs-unramified}
	We can compute the pairing for newforms in the case
	$\pi = \pi(\chi_1, \chi_2)$ with $\chi_1, \chi_2$ unramified. 
	Write $\alpha_i = \chi_i(\varpi)$ and $\beta_i = \wt{\chi_i}(\varpi)$.
	Then the newforms $W^\new$ and $\wt{W}^\new$ for $\pi$ and
	$\wt{\pi}$ take the values (as at the end of
	Section \ref{subsec:whittaker-kirillov}),
	\begin{align*}
		W^\new
			\paren{\begin{pmatrix} \varpi^n & 0 \\ 0 & 1 \end{pmatrix}}
			&= \frac{\alpha_1^{n + 1} - \alpha_2^{n + 1}}
				{\alpha_1 - \alpha_2} q^{\frac{-n}{2}}, \\
		\wt{W}^\new
			\paren{\begin{pmatrix} \varpi^n & 0 \\ 0 & 1 \end{pmatrix}}
			&= \frac{\beta_1^{n + 1} - \beta_2^{n + 1}}
				{\beta_1 - \beta_2} q^{\frac{-n}{2}}.
	\end{align*}
	Bringing this into the formula for $Z(s, W^\new, \wt{W}^\new)$
	(Equation~\ref{eq:zeta-integral-rs}) yields,
	\begin{align*}
		Z\paren{W^\new, \wt{W}^\new, s}
			= \, &\int_{N(F)Z(F) \bs \GL_2(F)}
				W^\new(g) \wt{W}^\new (\epsilon' g) f(g, s) dg \\
			= \, &\sum_n q^n W^\new
				\begin{pmatrix} \varpi^n & 0 \\ 0 & 1 \end{pmatrix}
				\wt{W}^\new
				\begin{pmatrix} \varpi^n & 0 \\ 0 & 1 \end{pmatrix} q^{-ns} \\
			= \, &\sum_{n = 0}^\infty
				\frac{\alpha_1^{n + 1} - \alpha_2^{n + 1}}{\alpha_1-\alpha_2}
				\frac{\beta_1^{n + 1} - \beta_2^{n + 1}}{\beta_1 - \beta_2}
				q^{-ns} \\
			= \, &\frac{1}{\paren{\alpha_1 - \alpha_2}\paren{\beta_1 - \beta_2}}
					\cdot \\
				\, &\paren{\frac{\alpha_1 \beta_1 q^{-s}}
					{1 - \alpha_1 \beta_1 q^{-s}}
					- \frac{\alpha_1 \beta_2 q^{-s}}{1 - \alpha_1 \beta_2 q^{-s}}
					- \frac{\alpha_2 \beta_1 q^{-s}}{1 - \alpha_2 \beta_1 q^{-s}}
					+ \frac{\alpha_2 \beta_2 q^{-s}}{1 - \alpha_2 \beta_2 q^{-s}}} \\
			= \, &\frac{1 - \alpha_1 \alpha_2 \beta_1 \beta_2 q^{-2s}}
					{\paren{1 - \alpha_1 \beta_1 q^{-s}}
					\paren{1 - \alpha_1 \beta_2 q^{-s}}
					\paren{1 - \alpha_2 \beta_1 q^{-s}}
					\paren{1 - \alpha_2 \beta_2 q^{-s}}} \\
			= \, &\frac{1 + q^{-s}}
				{\paren{1 - q^{-s}}\paren{1 - \frac{\alpha_1}{\alpha_2} q^{-s}}
				\paren{1 - \frac{\alpha_2}{\alpha_1} q^{-s}}} \\
			= \, &\paren{1 + q^{-s}} L\paren{\Ad(\pi), s}.
	\end{align*}

	Therefore, $\Psi(W^\new, \wt{W}^\new, s) = 1$
	by Equation~\ref{eq:normalized-zeta-rs}
	and
	\[
		\calp_\rs\paren{W^\new \otimes \wt{W}^\new}
			= Z\paren{W^\new, \wt{W}^\new, 1}
			= \paren{1 + q^{-1}} L\big(\Ad(\pi), 1\big).
	\]
\end{example}

\subsection{Theta liftings}
\label{sec:local-theta-liftings}

For a general orthogonal quadratic space
$(V, Q)$ of even dimension $m$ with the
bilinear form,
\[
	\pair{x, y} = Q(x + y) - Q(x) - Q(y),
\]
let $\go(V)$ be the group of similitudes on $V$ with norm map
$\nu: \go(V) \lra \G_m$. 
Let $G = \GL_2 \times_{\G_m} \go (V)$ be the fiber product of
$\nu$ and $\det: \GL_2 \lra \G_m$.
We consider $\SL_2$ and $\rmo(V)$ to be normal subgroups of
$G$ with respective quotients isomorphic to $\go(V)$ and 
$\GL_2(F)^+$,
the subgroup of elements $g \in \GL_2(F)$ such that
$\det g \in \nu(\go(V))$.
Let $\cals(V)$ be the space of Schwartz functions on $V$.

There is a Weil representation $r$ of $G(F)$ on $\cals(V)$ with
respect to $\psi: F \lra \C^\times$
(cf. \cite[Section~1]{waldspurger}, \cite[Section~2.1]{yzz}).
To describe this representation, we need the following
special elements in $\GL_2$:
\begin{align*}
	d(\lambda) &:= \begin{pmatrix} 1 & \\ & \quad \lambda \end{pmatrix}, \qquad
	a(\lambda) := \begin{pmatrix} \lambda & \\ & 1\end{pmatrix}, \\
	m(\lambda) &:= \begin{pmatrix} \lambda & \\ & \lambda^{-1} \end{pmatrix}, \qquad
	n(b) := \begin{pmatrix} 1 & b \\ & 1 \end{pmatrix}, \\
	w &:= \begin{pmatrix} & \, 1 \\ -1 & \end{pmatrix}.
\end{align*}
$G$ is generated by the elements $m(\lambda)$, $n(b)$, $w$, and
$(d(\nu(h)), h)$ for
$h \in \go(V)$.
The Weil representation $r$ is described as follows.
\begin{enumerate}
	\item For any $h \in \go(V)$ and
			$\Phi \in \cals(V)$,
			\[
				r\big(d\paren{\nu(h)},h\big)\cdot\Phi(x)
					=
				\verts{\nu(h)}^{-\frac{m}{4}}
					\Phi\paren{h^{-1}x}.
			\]
	\item For any $\lambda \in F^\times$,
		\[
			r\big(m(\lambda)\big) \cdot \Phi(x) =
				\eta_V(\lambda) \verts{\lambda}^{\frac{m}{2}} \Phi(\lambda x),
		\]
	where $\eta_V(\lambda) = (\lambda, (-1)^{m/2} \det(V))$.
	\item For any $b \in F$,
		\[
			r\big(n(b)\big) \cdot \Phi(x) = \psi\big(bQ(x)\big) \Phi (x).
		\]
	\item For $w$ as above,
		\[
			r(w) \cdot \Phi(x) = \gamma \cdot \wh \Phi(x),
		\]
		where $\gamma$ is an $8$-th root of unity
		and $\wh{\Phi}$ is the
		Fourier transform,
		\[
			\wh{\Phi}(x) = \int_{V} \Phi(y) \psi\paren{\pair{x, y}} dy.
		\]
\end{enumerate}

\begin{remark}
	\label{rem:d-lambda}
	We follow the convention of using $d(\nu(h))$
	instead of $a(\nu(h))$
	(following Harris--Kudla~\cite[Section~3.2]{harris-kudla-1991},
	\cite[Equation~1.1]{harris-kudla-2004},
	Yuan--Zhang--Zhang~\cite[Section~2.1.3]{yzz}
	in contrast to
	Jacquet--Langlands~\cite[Chapter~1]{jacquet-langlands}).
	Similar to Remark~\ref{rem:epsilon}, this is a notational decision
	that affects the appearance of $\chi(-1)$ in later calculations.
	We will also omit the $r$ where the context is clear,
	for example simply writing $w\Phi$ for $r(w) \cdot \Phi$
	and $n(b) \Phi$ for $r(n(b)) \cdot \Phi$.
\end{remark}

\begin{remark}
	\label{rem:weil-index}
	The $8$-th root of unity $\gamma$ in the action of
	$r(w)$ on a Schwartz function $\Phi$ is called the
	Weil index. It is dependent on $(V, Q)$, but
	is equal to $-1$ for nonsplit quaternion algebras
	and $1$ for split quaternion algebras
	(cf. \cite[Chapter~II]{weil})
	so it can be omitted for most of
	our purposes.
\end{remark}

\subsection{Whittaker functions of newforms and optimal forms}
\label{sec:quadratic-cases}
We start with the general quadratic space $V = (Ee, Q)$ with an
action of $E$, where $E$ is a semisimple algebra over $F$
and $Q$ is a multiple of the norm $\rmn = \rmn_{E/F}$ of $E$ over $F$.
Then $\go(V) = \pair{E^\times, \iota}$ where $\iota$ is an involution.
In this case, $\nu$ is the usual norm $\rmn$ of $E$ over $F$. 
 
Let $c_{E/F}$ denote the nontrivial involution of
$E/F$, and write
$\overline{x} := c_{E/F}(x)$ for $x \in E$.
Let $\chi: E^\times \lra \C^\times$ be a character.
For each $\Phi \in \cals(E)$, we obtain a Whittaker
function supported by the subgroup
$\GL_2(F)^+$ of matrices with determinant in
$\rmn(E^\times)$.
Thus $\GL_2(F)^+ = \GL_2(F)$ if
$E = F \oplus F$; otherwise, $\GL_2(F)^+$ is an
index-$2$ subgroup of $\GL_2(F)$.
More precisely, for $g \in \GL_2(F)^+$, write
$g = d(Q(t_0e)^{-1})g_1$ with
$t_0 \in E^\times$ and $g_1 \in \SL_2(F)$.
Then
\begin{equation}
	\label{eq:whittaker}
	W(g,\chi,\Phi)
		=
	\verts{\det g}^{-\frac{1}{2}}
	\int_{E^1}
		r\paren{g_1}\Phi\paren{tt_0e}
		\chi^{c_{E/F}}
			\paren{t_0^{-1}t^{-1}}dt.
\end{equation}
Here $E^1$ is the subgroup of $E^\times$ with norm
$1$, and $dt$ is the Haar measure on $E^1$
normalized by $\vol(\calo_E^1) = 1$.
The corresponding Kirillov functions are given by
\begin{equation}
	\label{eq:kappa-chi}
	\kappa(x,\chi,\Phi)
		=
	\verts{x}^{\frac{1}{2}}
	\int_{E^1}
		\Phi\paren{t_0te}\chi\paren{tt_0}dt,
\end{equation}
where $x = Q(t_0e)$.
\begin{remark}
	Note that the above construction is compatible with the
	construction in the global situation,
	\cite[Equations 5.2 and 5.3]{zhang-hv}, by the decomposition
	given by \cite[Equations 5.4 and 5.5]{zhang-hv}.
\end{remark}

The normalization in Equation~\ref{eq:whittaker}
gives
$W(n(b)g, \chi, \Phi)
= \psi(b)W(g, \chi, \Phi)$.
The functions $W(g, \chi, \Phi)$ form the explicit
local theta lifting $\theta(\chi, \psi)$ on
$\GL_2(F)^+$. The map
\[
	\begin{tikzcd}[row sep = tiny]
		\cals(V) \arrow[r]
			& \theta(\chi, \psi) \\
		\Phi \arrow[r, mapsto]
			& W(\,\cdot\,, \chi, \Phi)
	\end{tikzcd}
\]
is equivariant for
$E^\times \times_{\G_m}\GL_2(F)$. Equivalently, it
defines an element of
\[
	\Hom_{E^\times \times_{\G_m}\GL_2(F)}
		\paren{
			\cals(V) \otimes \chi,
			\theta(\chi, \psi)}.
\]

The space $\theta(\chi, \psi)$ is stable under right
translation by $\GL_2(F)^+$, and we regard it as a
space of functions on $\GL_2(F)$ supported on
$\GL_2(F)^+$.
If $E = F \oplus F$, set
$\calw(\chi, \psi) := \theta(\chi, \psi)$.
Otherwise, choose
$h \in \GL_2(F) - \GL_2(F)^+$, and let
$\calw(\chi, \psi)$ be the space of functions
\[
	W(g) = W_1(g) + W_2(gh),
\]
where $W_1, W_2 \in \theta(\chi, \psi)$.
In either case, $\calw(\chi, \psi)$ is an
irreducible representation of $\GL_2(F)$,
denoted by $\pi(\chi)$.
This space has the following properties.
\begin{enumerate}
	\item The central character of $\pi(\chi)$ is
		$\omega := \eta \cdot \restr{\chi}{F^\times}$;
	\item If $E = F \oplus F$ and $\chi = (\chi_1, \chi_2)$ then
		$\pi(\chi) = \pi(\chi_1, \chi_2)$ is a principal series;
	\item If $\chi = \omega \circ \rmn$, then
		$\chi = (\omega, \omega \cdot \eta)$ is a principal series;
	\item If $E$ is not split and $\chi$ does not factor through the
		norm map, then $\pi(\chi)$ is supercuspidal.
\end{enumerate}

\subsubsection*{Newforms}
We assume that $V = (E, \rmn)$ with $E$ a
semisimple algebra over $F$.
For the induced representation $\pi(\chi)$, the
newform $W^\new_\chi$ can be constructed explicitly
from the Whittaker function $W(g, \chi, \Phi_\chi)$ by the following
two steps.
\begin{enumerate}
	\item Define $\Phi_\chi$ as in \cite[Section~5.2]{zhang-hv}.
		\begin{enumerate}
			\item If $E$ is a field extension and $\chi$ is unramified,
				then define $\Phi_\chi$ to be the characteristic function of
				$\calo_E$;
			\item If $E$ is a field extension and $\chi$ is ramified, then
				define $\Phi_\chi$ to be the restriction of $\chi^{-1}$ on
				$\calo_E^\times$;
			\item If $E = F \oplus F$,
				$\chi = (\chi_1, \chi_2)$, then define
				$\Phi_\chi = \Phi_{\chi_1} \otimes \Phi_{\chi_2}$,
				where for each $i$,
				\[
					\Phi_{\chi_i} = \begin{cases}
						\restr{\one}{\calo_F} & \text{if $\chi_i$ is unramified,} \\
						\restr{\chi_i^{-1}}{\calo_F^\times} & \text{otherwise.}
					\end{cases}
				\]
		\end{enumerate}
		Then $W(g, \chi, \Phi_\chi)$ is already a newform $W^\new_\chi(g)$
		when $E/F$ is not ramified.
	\item If $E/F$ is ramified, then the newform in $\pi(\chi)$ has the
		form,
		\[
			W^\new_\chi(g) = W\paren{g, \chi, \Phi_\chi}
				+ W\paren{ga(\epsilon), \chi, \Phi_\chi},
		\]
		where $\epsilon\in \calo_F^\times \setminus \rmn(\calo_{E}^\times)$.
\end{enumerate}

This shows, in particular, the following formula: 
\begin{prop}
	\label{prop:theta-new}
	The inert product of newforms can be computed from theta forms as
	follows,
	\[
		\calp_{\rs} \paren{W^\new_\chi\otimes W^\new _{\chi^{-1}}}
			= e_{E/F} \calp_{\rs}
				\paren{W\paren{g, \chi^{-1}, \Phi_{\chi^{-1}}} \otimes
				W\paren{g, \chi, \Phi_\chi}},
	\]
	\noindent
	where 
	\[
		e_{E/F} := \begin{cases}
			1 & \text{if $E/F$ is unramified,} \\
			2 & \text{if $E/F$ is ramified}.
		\end{cases}
	\]
\end{prop}

\subsubsection*{Optimal forms}
We now construct a distinguished element
$W^\opt \in \calw(\chi, \psi)
\otimes \calw(\chi^{-1}, \psi)$,
which is the local Whittaker analogue of the optimal
form defined in \cite[Definition~6.1]{zhang-hv}.
The adjective ``optimal'' refers to optimal embeddings.
An embedding $\iota: E \longhookrightarrow B$ is optimal
with respect to orders $\calo_c \subset E$ and
$\calo_B \subset B$ if
$\iota(E) \cap \calo_B = \iota(\calo_c)$;
see \cite[Remark~6.2]{zhang-hv}.
Thus the order cut out inside $E$ is exactly
$\calo_c$, rather than a larger order.
As Lemma~\ref{lem:opt-pm} below shows, the resulting
two-variable vector comes from the characteristic
function of $\End_{\calo_F}(\calo_c)$.
This realization is what makes its Rankin--Selberg
period amenable to explicit calculation.

Let $\xi := \chi^{1 - c_{E/F}}$ be the antinorm
character, so
$\xi(x) = \chi(x / \overline{x})$.
Under the isomorphism
\[
	\begin{tikzcd}[row sep = tiny]
		E^\times / F^\times
			\arrow[r, "\widesim{}"]
			& E^1 \\
		x \arrow[r, mapsto]
			& \overline{x} / x
	\end{tikzcd}
\]
the character $\xi$ corresponds to
$\restr{\chi^{-1}}{E^1}$.
Then $\xi$ is a ring class character.

Let $\varpi$ be a uniformizer of $F$.
The least non-negative integer $o(\xi)$ such that
$\xi$ is trivial on
$(\calo_F + \varpi^{o(\xi)}\calo_E)^\times$
is called the order of $\xi$.
Set $c := c(\xi) := \varpi^{o(\xi)}$, and write
\[
	\calo_c := \calo_F + c\calo_E
\]
for the associated order of $E$.

Let $\delta \in \calo_c$ be a generator of the different ideal $\cald$ of
$\calo_c$, namely the ideal generated by $x - \overline{x}$ for all
$x \in \calo_c$. Then for each $a \in \calo_c / \delta$, we define the
function $\Phi_a^\opt$ to be the characteristic function of 
\[
	\calo_c + \frac{a}{\delta} \subset E.
\]
Letting $\epsilon :=
\begin{psmallmatrix} -1 & 0 \\ 0 & 1\end{psmallmatrix}$,
the optimal form is defined as follows.

\begin{defn}
\label{def:opt-function}
	We define the one-variable optimal form $W_a^\opt$
	on $\GL_2(F)$,
	for $a \in \calo_c/\delta$, and the
	two-variable optimal form $W^\opt$ on $\GL_2(F) \times \GL_2(F)$,
	\begin{align*}
		W_a^\opt(g) &:= W\paren{g, \chi, \Phi^\opt_a}, \\
		W^\opt\paren{g_1, g_2} &:= \sum_{a \in \calo_c/\delta}
			W_a^\opt\paren{g_1} \otimes W_{-a}^\opt\paren{g_2\epsilon}.
	\end{align*}
\end{defn}

\subsubsection*{Comparison of models}
Let $V = (Ee, Q)$ and $V' = (Ee', Q')$ be
quadratic spaces with an action of $E$, and let
$\iota: V' \iso V$ be the $E$-linear isomorphism
satisfying $\iota(e') = e$.
For $\Phi \in \cals(V)$, set
$\iota^*\Phi := \Phi \circ \iota \in \cals(V')$.
Define
\[
	Q(\iota) := \frac{Q(e)}{Q'(e')} \in F^\times.
\]

Write $\kappa_V$ and $\kappa_{V'}$ for the
corresponding Kirillov functions.
If $x = Q(t_0e)$, then
$Q'(t_0e') = xQ(\iota)^{-1}$ and
$\iota^*\Phi(tt_0e') = \Phi(tt_0e)$.
Equation~\ref{eq:kappa-chi} therefore gives
\begin{align*}
	&\kappa_{V'}\paren{
		xQ(\iota)^{-1},
		\chi^{c_{E/F}},
		\iota^*\Phi} \\
	&\qquad
		= \verts{xQ(\iota)^{-1}}^{\frac{1}{2}}
			\int_{E^1}
				\Phi\paren{tt_0e}
				\chi\paren{tt_0}dt \\
	&\qquad
		= \verts{Q(\iota)}^{-\frac{1}{2}}
			\kappa_V\paren{
				x,\chi^{c_{E/F}},\Phi}.
\end{align*}
Thus
\begin{equation}
	\label{eq:kk}
	\kappa_V\paren{x,\chi^{c_{E/F}},\Phi}
		= \verts{Q(\iota)}^{\frac{1}{2}}
			\kappa_{V'}\paren{
				xQ(\iota)^{-1},
				\chi^{c_{E/F}},
				\iota^*\Phi}.
\end{equation}

For example, if we compare the Whittaker functions defined by two
opposite spaces $V_{\pm} := (E, \pm \rmn)$, we get two Whittaker
functions $W_\pm(g, \chi^{c_{E/F}}, \Phi)$ for $\Phi \in \cals(E)$.
We use the identity map $\iota: V_+ \lra V_- $ for the quadratic
space, so $Q(\iota) = -1$.
Then
\[
	W_-\paren{g, \chi^{c_{E/F}}, \Phi} = W_{+}\paren{g\epsilon, \chi^{c_{E/F}}, \Phi},
\]
where $\epsilon = \begin{pmatrix} -1 & 0 \\ 0 & 1 \end{pmatrix}$.

Then from Definition~\ref{def:opt-function}, we can instead write 
\begin{equation}
	\label{eq:opt-pm}
	W^\opt\paren{g_1, g_2}
		= \sum_{a\in \calo_c / \delta} W_+\paren{g_1, \chi, \Phi_a^\opt}
			\cdot W_-\paren{g_2, \chi^{-1}, \Phi_{-a}^\opt}.
\end{equation}

We identify $B = \End_F(E)$ as usual and take an optimal order
$\calo_B^\opt = \End_{\calo_F}(\calo_c)$. 
Let $\Phi^\opt$ be the characteristic function of $\calo_B^\opt$.

\begin{lem}
	\label{lem:opt-pm}
	We have the following identity in
	$\cals(B) = \cals(V_+) \otimes \cals(V_-)$,
	\[
		\Phi^\opt = \sum_a \Phi^\opt_a \otimes \Phi^\opt_{-a}.
	\]
\end{lem}

\begin{proof}
	Choose $\tau \in \calo_c$ such that
	$\calo_c = \calo_F + \calo_F\tau$.
	After replacing $\delta$ by an associate, we may
	assume that $\delta = \tau - \overline{\tau}$.
	This only relabels the classes in $\calo_c/\delta$.

	Write $B = E \oplus Ej$, where
	$jz = \overline{z}j$ for every $z \in E$.
	Then $x + yj$ acts on $E$ by
	$z \mapsto xz + y\overline{z}$.
	It belongs to
	$\calo_B^\opt = \End_{\calo_F}(\calo_c)$
	if and only if its values on the basis $1, \tau$
	belong to $\calo_c$, i.e.
	\begin{align*}
		s &:= x + y \in \calo_c, \\
		t &:= x\tau + y\overline{\tau}
			\in \calo_c.
	\end{align*}
	Solving for $x$ and $y$ gives
	\begin{align*}
		x &= \frac{t - s\overline{\tau}}{\delta}, \\
		y &= \frac{s\tau - t}{\delta}.
	\end{align*}
	Hence $x, y \in \delta^{-1}\calo_c$ and
	$x + y \in \calo_c$.
	Let $a$ be the class of $\delta x$ in
	$\calo_c/\delta$.
	Then $x \in \calo_c + a/\delta$ and
	$y \in \calo_c - a/\delta$.

	Conversely, if
	$x \in \calo_c + a/\delta$ and
	$y \in \calo_c - a/\delta$, then
	$x + y \in \calo_c$ and
	$x\tau + y\overline{\tau} \in \calo_c$.
	Thus
	\[
		\calo_B^\opt
			=
		\bigsqcup_{a \in \calo_c/\delta}
		\paren{\calo_c + \frac{a}{\delta}}
		+
		\paren{\calo_c - \frac{a}{\delta}}j.
	\]
	The class $a$ is determined by $x$ modulo
	$\calo_c$, so the union is disjoint.
	Taking characteristic functions gives the
	claimed identity.
\end{proof}


\subsection{Optimal forms for locally dihedral modular forms}
\label{sec:locally-dihedral}
We now extend the definition of optimal forms
given in \cite[Section~6]{zhang-hv} from dihedral 
weight-$1$ modular forms
to locally dihedral weight-$1$ modular forms.

\subsubsection{Locally dihedral forms}
Let $f$ be a newform of weight $1$ and level $N$,
and let $\pi_f$ be its automorphic representation.
We call $f$ \emph{locally dihedral} if, for every
prime $p \mid N$, there are a quadratic semisimple
$\Q_p$-algebra $K_p$ and a character
$\chi_p:K_p^\times \lra \C^\times$ such that
\[
	\pi_{f,p} \cong \pi(\chi_p).
\]
If $K_p/\Q_p$ is a field extension, this is
equivalent under local Langlands to
\[
	\rho_{f,p}
		\cong
	\Ind_{G_{K_p}}^{G_{\Q_p}}\paren{\chi_p}.
\]
If $K_p = \Q_p \oplus \Q_p$, write
$\chi_p = (\chi_{p,1},\chi_{p,2})$; then
$\pi(\chi_p)$ is the corresponding principal series.

This condition is automatic at odd primes,
hence $f$ is locally dihedral if $N$ is odd.
If $2 \mid N$, the condition at $2$ also holds
when $\rho_f$ is $2$-ordinary.
Set
\begin{align*}
	\Q_N := \prod_{p \mid N}\Q_p, \qquad
	K_N := \prod_{p \mid N}K_p, \qquad
	\chi_N := \prod_{p \mid N}\chi_p.
\end{align*}
Then $K_N$ is a quadratic semisimple algebra
over $\Q_N$.
These objects play the same role for locally
dihedral forms as $K/\Q$ and $\chi$
do for dihedral forms,
so we will often denote them without the subscript $N$.
  
\subsubsection{Construction of optimal forms}
Let $f$ be a locally dihedral form of weight $1$.
We define a two-variable
optimal form $f^\opt(z_1, z_2)$
in the space of $f(z_1)f^*(z_2)$ generated by
$\GL_2(\Q_N)$, analogously to \cite[Equation 6.1]{zhang-hv} but using
Whittaker functions.
More precisely, let $\varphi$ and $\varphi^*$ be the automorphic
avatars of $f$ and $f^*$, and let $W(g)$ and $W^*(g)$ be their
Whittaker coefficients:
\begin{align*}
	\varphi(g) &= \sum_{a \in \Q^\times}
		W \paren{\begin{pmatrix}a & \\ & 1 \end{pmatrix} g}, \\ 
	\varphi^*(g) &= \sum_{a \in \Q^\times}
		W^* \paren{\begin{pmatrix}a & \\ & 1 \end{pmatrix} g}.
\end{align*}

Then we have decompositions into products of local newforms:
\begin{align*}
	W(g) &= \prod_{p \leq \infty} W_p\paren{g_p}, \\
	W^*(g) &= \prod_{p \leq \infty} W_p^*\paren{g_p}.
\end{align*}
  
To construct $f^\opt(z_1, z_2)$, it suffices to construct the
Whittaker coefficients $W^\opt(g_1, g_2)$
of its automorphic avatar $\varphi^\opt (g_1, g_2)$: 
\[
	\varphi^\opt\paren{g_1, g_2} = \sum_{a, b \in \Q^\times}
		W^\opt \paren{\begin{pmatrix}a & \\ & 1 \end{pmatrix} g_1,
			\begin{pmatrix}b & \\ & 1 \end{pmatrix} g_2}.
\]
We will construct the local Whittaker functions $W^\opt_p$ and then
put them together:
\[
	W^\opt\paren{g_1, g_2} := \prod_{p \leq \infty} W^\opt_p\paren{g_{1, p}, g_{2, p}}.
\]
For $p \nmid N$, set
\[
	W^\opt_p\paren{g_1, g_2} := W_p\paren{g_1} W^*_p\paren{g_2}.
\]
For $p \mid N$, we want to construct an optimal
element in
$\calw(\chi_p,\psi_p)
\otimes\calw(\chi_p^{-1},\psi_p)$.
We specialize the preceding local construction to
$F = \Q_p$, $E = K_p$, $\chi = \chi_p$, and
$\psi = \psi_p$.

Set $\xi_p := \chi_p^{1 - c_{E/F}}$, so
$\xi_p(x) = \chi_p(x / \overline{x})$.
Under the isomorphism
\[
	\begin{tikzcd}[row sep = tiny]
		K_p^\times / \Q_p^\times
			\arrow[r, "\widesim{}"]
			& K_p^1 \\
		x \arrow[r, mapsto]
			& \overline{x} / x
	\end{tikzcd}
\]
the character $\xi_p$ corresponds to
$\restr{\chi_p^{-1}}{K_p^1}$.

Let $\varpi_p$ be a uniformizer of $\Q_p$.
Let $o(\xi_p)$ be the least non-negative integer
such that $\xi_p$ is trivial on
\[
	\paren{
		\Z_p + \varpi_p^{o(\xi_p)}
			\calo_{K,p}}^\times.
\]
Set $c_p := \varpi_p^{o(\xi_p)}$ and
$\calo_{c,p} := \Z_p + c_p\calo_{K,p}$.
Let $\delta_p \in \calo_{c,p}$ generate the different
ideal $\cald_p$ of $\calo_{c,p}$.
For each $a \in \calo_{c,p} / \delta_p$, let
$\Phi_{a,p}^\opt$ be the characteristic function of
\[
	\calo_{c,p} + \frac{a}{\delta_p}
		\subset K_p.
\]

We then define the one-variable optimal function
$W_{a,p}^\opt$ and the two-variable optimal function
$W_p^\opt$ by
\begin{align*}
	W_{a,p}^\opt(g)
		&:= W\paren{g, \chi_p, \Phi_{a,p}^\opt}, \\
	W_p^\opt\paren{g_1, g_2}
		&:= \sum_{a \in \calo_{c,p} / \delta_p}
			W_{a,p}^\opt\paren{g_1}
			\otimes
			W_{-a,p}^\opt\paren{g_2\epsilon}.
\end{align*}

\begin{rem}
	The purpose of $\varphi^\opt$ is to replace the
	newform tensor by a test vector whose periods are
	explicitly computable. The non-vanishing result below
	is the local input showing that the canonical period
	of the global optimal form is nonzero.
	The comparison theorem in \cite[Theorem~6]{zhang-hv} transfers a
	Harris--Venkatesh identity for $\varphi^\opt$ to the
	newform tensor after multiplication by fixed integral
	factors.
	To establish the Harris--Venkatesh and Stark
	conjectures compatibly, as predicted in
	\cite[Conjecture~6]{zhang-hvs},
	the natural intermediate statement is the
	existence of a unique element
	$u_{\varphi^\opt}
	\in \calu\paren{\Ad(\rho_f)} \otimes \Q$
	satisfying
	$\calp_\rs\paren{\varphi^\opt}
	= \Reg_\R\paren{u_{\varphi^\opt}}$ and
	$\calp_\hv\paren{\varphi^\opt}
	= \Reg_{\F_p^\times}
		\paren{u_{\varphi^\opt}}$
	for all sufficiently large $p$.
\end{rem}


\section{Non-vanishing and rationality of local Rankin--Selberg periods}
\label{sec:non-vanishing}

\subsection{Non-vanishing on unitary representations}
Assume that $\pi$ is unitary, so
$\wt{\pi} \cong \overline{\pi}$.
For $W \in \calw(\pi, \psi)$, define
$W^\epsilon \in \calw(\overline{\pi}, \psi)$ by
\begin{equation}
	\label{eq:bar-W}
	W^\epsilon(g) := \overline{W(\epsilon g)}.
\end{equation}
The conjugation compensates for the change from $\psi$
to $\psi^{-1}$ under left translation by $\epsilon$.

\begin{prop}
	\label{prop:R-unitary}
	If $\pi$ is unitary and $W \neq 0$, then
	$\calp_{\rs}\paren{W \otimes W^\epsilon}$
	is nonzero and real.
\end{prop}

\begin{proof}
	For $x \in F$, the identity
	$\epsilon n(x) = n(-x)\epsilon$ gives
	\begin{align*}
		W^\epsilon\paren{n(x)g}
			&= \overline{W\paren{\epsilon n(x)g}} \\
			&= \overline{W\paren{n(-x)\epsilon g}} \\
			&= \overline{\psi(-x)}W^\epsilon(g) \\
			&= \psi(x)W^\epsilon(g).
	\end{align*}
	Thus $W^\epsilon$ belongs to
	$\calw(\overline{\pi}, \psi)$, as asserted.

	Let $\omega$ be the central character of $\pi$.
	Since $\epsilon\epsilon' = -\Id$, one has
	\[
		W^\epsilon(\epsilon'g)
			= \overline{\omega(-1)}
				\overline{W(g)}.
	\]
	Equation~\ref{eq:zeta-integral-rs} therefore gives
	\[
		\calp_{\rs}\paren{W \otimes W^\epsilon}
			= \overline{\omega(-1)}
				\int_{N(F)Z(F) \bs \GL_2(F)}
				\verts{W(g)}^2 f(g, 1) dg.
	\]
	The integral converges absolutely, and its integrand
	is non-negative and not identically zero.
	Finally, $\omega(-1) \in \set{\pm 1}$, so the
	period is nonzero and real.
\end{proof}

\subsection{Non-vanishing on optimal forms}

One of our main objectives is to study the period
$\calp_{\rs}(W^\opt)$.
First, we prove a non-vanishing result.
\begin{prop}
	\label{prop:opt-unitary}
	If $\chi$ is unitary, then 
	\[
		W^\opt\paren{g_1, g_2} = \chi(-1) \sum_{a \in \calo_c / \delta}
			W_a^\opt\paren{g_1} \otimes W_{a}^{\opt, \epsilon}\paren{g_2}.
	\]
	Furthermore, $W_1^\opt \neq 0$ and
	$\calp_{\rs}(W^\opt) \in \R^\times$.
\end{prop}

\begin{proof}
	First, let us write the Kirillov functions for $W^\opt$:
	\[
		\kappa^\opt\paren{x_1, x_2} :=
			\sum_{a \in \calo_c / \delta} \kappa\paren{x_1, \chi, \Phi_a^\opt}
			\cdot \kappa\paren{-x_2, \chi^{-1}, \Phi_{-a}^\opt}.
	\]
	Since
	$\Phi_{-a}^\opt(x)=\Phi_a^\opt(-x)$, we can use
	Equation~\ref{eq:kappa-chi} and change variables by
	$t\mapsto-t$ to obtain
	\[
		\kappa^\opt\paren{x_1,x_2}
			=
		\chi(-1)\sum_{a\in\calo_c/\delta}
			\kappa\paren{x_1,\chi,\Phi_a^\opt}
			\cdot
			\kappa\paren{-x_2,\chi^{-1},\Phi_a^\opt}.
	\]
	Since $\chi$ is unitary,
	\[
		\kappa\paren{-x_2,\chi^{-1},\Phi_a^\opt}
			=
		\overline{
			\kappa\paren{-x_2,\chi,\Phi_a^\opt}}
			=
		\kappa^\epsilon
			\paren{x_2,\chi,\Phi_a^\opt},
	\]
	where
	$\kappa^\epsilon(x,\chi,\Phi_a^\opt)$ is the
	Kirillov function associated to
	$W^\epsilon(g,\chi,\Phi_a^\opt)$ in
	Equation~\ref{eq:bar-W}.
	Thus,
	\[
		\kappa^\opt\paren{x_1, x_2} :=
			\chi(-1) \sum_{a \in \calo_c / \delta}
			\kappa\paren{x_1, \chi, \Phi_a^\opt} \cdot
			\kappa^\epsilon\paren{x_2, \chi, \Phi_{a}^\opt}.
	\]
	It follows that
	\[
		\calp_{\rs}\paren{W^\opt}
			= \chi(-1) \sum_{a \in \calo_c/\delta}
			\calp_{\rs}\paren{W_a^\opt \otimes W_a^{\opt, \epsilon}}.
	\]
	For the non-vanishing of $W_1^\opt$,
	take $x = \rmn(\delta^{-1})$ in
	Equation~\ref{eq:kappa-chi} to obtain 
	\[
		\kappa_1^\opt(x)
			= \chi(\delta)^{-1} \verts{\delta}^{-\frac{1}{2}}
				\int_{\paren{1 + \delta \calo_c}^1} \xi(t)
			= \chi(\delta)^{-1} \verts{\delta}^{-\frac{1}{2}}
			\vol\paren{\paren{1 + \delta \calo_c}^1} \neq 0.
	\]
	
	The last part of the proposition
	follows from the previous two parts and
	Proposition \ref{prop:R-unitary}.
\end{proof}

\subsection{Rationality}

Let $E/F$ be a semisimple algebra of degree $2$ and
$\chi$ be a character of $E^\times$.
We have constructed Whittaker models $\calw(\chi^\pm)$ for
$\GL_2(F)^+$ via theta liftings on quadratic spaces
$V_1 = (E, \rmn)$ and $V_2 = (E, -\rmn)$.
More precisely, for functions
$\Phi_i \in \cals(V_i)$,
Equation~\ref{eq:whittaker} gives
\begin{align*}
	W_1(g)
		&= W_+\paren{g,\chi,\Phi_1} \\
		&= \verts{\det g}^{-\frac{1}{2}}
			\int_{E^1}
				r\paren{g_1}
				\Phi_1\paren{t_0^{-1}t_1^{-1}}
				\chi^{c_{E/F}}\paren{t_0t_1}dt_1, \\
	W_2(\epsilon'g)
		&= W_-\paren{\epsilon'g,\chi^{-1},\Phi_2} \\
		&= \verts{\det g}^{-\frac{1}{2}}
			\int_{E^1}
				r\paren{g_1}
				\Phi_2\paren{t_0^{-1}t_2^{-1}}
				(\chi^{-1})^{c_{E/F}}
					\paren{t_0t_2}dt_2.
\end{align*}
Here $t_0 \in E^\times$ and
$g_1 \in \SL_2(F)$
satisfy
$g = d\paren{\rmn(t_0)}g_1$ and
$\epsilon'g = d\paren{-\rmn(t_0)}g_1$.

With the above, we compute the zeta integral from
Equation~\ref{eq:zeta-integral-rs}
to obtain
\begin{align*}
	Z(W_1, W_2, s)
		= \, &\int_{N(F)Z(F) \bs \GL_2(F)^+} W_1(g) W_2(\epsilon' g) f(g, s)dg \\
		= \, &\int_{N(F)Z(F) \bs \GL_2(F)^+}\verts{\det g}^{-1} f(g, s) dg \\
			\, &\cdot \int_{E^1 \times E^1} r\paren{g_1}
			\Phi\paren{t_0^{-1} \paren{t_1^{-1} + t_2^{-1}j}}
			\chi^{c_{E/F}}\paren{\frac{t_1}{t_2}} dt_1 dt_2.
\end{align*}
Normalize the Haar measure $dk$ on
$\SL_2(\calo_F)$ by
$\vol\paren{\SL_2(\calo_F)} = 1$, and set
\[
	\wt{\Phi}(x)
		:= \int_{\SL_2(\calo_F)} r(k)\Phi(x) dk.
\]
Then the integral expression of $Z(W_1, W_2, s)$ becomes
\[
	\int_{N(F) Z(F) \bs \GL_2(F)^+ / \SL_2(\calo_F)}
		\verts{\det g}^{-1} f(g, s) dg
		\int_{E^1 \times E^1} r\paren{g_1}
		\wt{\Phi}\paren{t_0^{-1}\paren{t_1^{-1} + t_2^{-1}j}}
		\chi^{c_{E/F}}\paren{\frac{t_1}{t_2}} dt_1 dt_2.
\]
Using the Iwasawa decomposition
\[
	\GL_2(F)^+
		= N(F)Z(F)
			d\paren{\rmn(E^\times)}
			\SL_2(\calo_F),
\]
the first integral becomes an integral over
$\rmn(E^\times)$ with respect to the restriction
of the Haar measure on $F^\times$.

Keep the Haar measure on $E^1$ normalized by
$\vol(\calo_E^1) = 1$, and endow $E^\times$
with the Haar measure compatible with the exact
sequence
\[
	1 \lra E^1 \lra E^\times
		\xrightarrow{\rmn}
		\rmn(E^\times) \lra 1.
\]
Thus
\[
	\vol(\calo_E^\times)
		=
	\kappa_{E/F}
		:=
	\sbrac{\calo_F^\times:
		\rmn(\calo_E^\times)}^{-1}.
\]
In particular, $\kappa_{E/F} = 1$ if $E/F$ is
split or unramified, while
$\kappa_{E/F} = 1/2$ if $E/F$ is ramified.

Let
\[
	T
		:=
	E^\times \times_{F^\times} E^\times
		=
	\set{\paren{a,b} \in E^\times \times E^\times
		\Mid \rmn(a) = \rmn(b)}.
\]
Setting
$\sigma(a,b) := \chi^{c_{E/F}}(b/a)$, the unfolding allows us to
rewrite the expression for $Z(W_1, W_2, s)$ as
\[
	\int_T \verts{\nu(t)}^s
		\wt{\Phi}\paren{t(1 + j)}
		\sigma(t) dt.
\]
The isomorphism
\[
	\begin{tikzcd}[row sep = tiny]
		E^\times\times E^1 \arrow[r, "\widesim{}"]
			& T \\
		\paren{t_1,t_2} \arrow[r, mapsto]
			& \paren{t_1,t_1t_2}
	\end{tikzcd}
\]
identifies the measure arising from the unfolding
with the product of the above measures on
$E^\times$ and $E^1$.
For
$\Phi = \sum_i \Phi_{1,i} \otimes \Phi_{2,i}$ in
$\cals(B) = \cals(V_+) \otimes \cals(V_-)$, set
\[
	W_\Phi\paren{g_1, g_2}
		:= \sum_i
			W_+\paren{g_1, \chi, \Phi_{1,i}}
			\cdot W_-\paren{g_2, \chi^{-1}, \Phi_{2,i}}.
\]
This does not depend on the chosen expression for $\Phi$.
Recall that $\xi = \chi^{1 - c_{E/F}}$ can be viewed as the
restriction $\restr{\chi^{-1}}{E^1}$.

\begin{prop}
	\label{prop:rs-phi}
	For every $\Phi \in \cals(B)$,
	\[
		Z\paren{W_\Phi, s}
			= Z\paren{\Phi, s}
			:= \int_{E^\times} \verts{t_1}^s
				\int_{E^1}
					\wt{\Phi}\paren{t_1\paren{1 + t_2j}}
					\xi(t_2) dt_2 dt_1.
	\]
\end{prop}

\begin{proof}
	For a pure tensor $\Phi = \Phi_1 \otimes \Phi_2$,
	this is the unfolding calculated above.
	Both sides depend linearly on $\Phi$, so the formula
	follows for every $\Phi \in \cals(B)$.
\end{proof}

\begin{cor}
	\label{cor:rat}
	Let $\Q(\xi, \wt{\Phi})$ be the subfield of $\C$
	generated by the values of $\xi$ and $\wt{\Phi}$.
	Then $\calp_{\rs}(W_\Phi)
	\in \Q\paren{\xi, \wt{\Phi}}$.
\end{cor}

\begin{proof}
	The functions $\xi$ and $\wt{\Phi}$ are locally
	constant, and $\wt{\Phi}$ is compactly supported.
	The Haar measures of compact open subsets are rational
	under the normalizations above. Decomposing the integral
	in Proposition~\ref{prop:rs-phi} by valuations therefore
	expresses $Z(\Phi, s)$ as a rational function in $q^{-s}$
	with coefficients in $\Q(\xi, \wt{\Phi})$.
	It is regular at $s = 1$ by absolute convergence, so the
	result follows from
	$\calp_{\rs}(W_\Phi) = Z(W_\Phi, 1)$.
\end{proof}

We now turn to the zeta integral
for the optimal form $W^\opt$ defined by
Equation~\ref{eq:opt-pm}.
By Lemma \ref{lem:opt-pm}, $W^\opt$ is defined by 
an {\em optimal function} $\Phi^\opt$ with respect to $\xi$ which is
the characteristic function of,
\[
	\calo_B^\opt = \calo_c + \calo_c \frac{1 - j}{\varpi - \overline{\varpi}},
\]
where $\varpi$ is a uniformizer of $E$.
Note that $\wt{\Phi}^\opt = \Phi^\opt$, and that 
$t_1 + t_2j \in \calo_B$ if and only if 
\begin{align*}
	t_1 &\in \delta^{-1} \calo_c, \\
	t_2 &\in -1 + t_1^{-1} \calo_c.
\end{align*}
Thus we have the following expression.
\begin{prop}
	\label{prop:rs-opt}
	\[
		Z\paren{W^\opt, s}
			= Z\paren{\Phi^\opt, s}
			= \chi(-1) \int_{\delta^{-1}\calo_c} \verts{t_1}^s dt_1
				\int_{\paren{1 + t_1^{-1}\calo_c}^1} \xi\paren{t_2} dt_2.
	\]
\end{prop}

Combining Propositions \ref{prop:rs-opt} and \ref{prop:opt-unitary},
we obtain the following rationality statement.
\begin{cor}
	\label{cor:opt-rat}
	Let $\Q(\xi + \xi^{-1})$ denote the ring generated by values of
	$\xi + \xi^{-1}$. Then
	\[
		\calp_\rs\paren{W^\opt} \in \Q\paren{\xi + \xi^{-1}}.
	\]
	Furthermore, if $\xi$ is unitary, then
	\[
		\calp_\rs\paren{W^\opt} \in \Q\paren{\xi + \xi^{-1}}^\times.
	\]
\end{cor}

We thereby obtain the rationality
part of Theorem \ref{thm:ratio} by combining
Corollary \ref{cor:opt-rat} and the non-vanishing of
$\calp_\rs(W^\new)$
(since new vectors are generators).


\section{Local Rankin--Selberg periods of optimal forms}
\label{sec:inner-product-opt}

The local factor in Theorem~\ref{thm:ratio} is the
quotient of the Rankin--Selberg periods of the new and
optimal vectors. We first calculate the denominator
$\calp_\rs(W^\opt)$. Proposition~\ref{prop:rs-opt}
reduces this period to an integral determined by the
order $\calo_c$, its different $\delta$, and the
restriction of $\xi$ to $E^1$. The three cases below
separate the contributions from the ramification of
$E/F$ and of $\xi$.

We calculate $Z(W^\opt, s)$ and
$\calp_\rs(W^\opt)$ when $\xi$ is respectively
unramified, ramified with $E/F$ inert, and ramified
with $E/F$ split.

\subsection{Unramified calculation}

In this subsection, we calculate $\calp_\rs(W^\opt)$ when
$\xi$ is unramified (i.e. as a character from $E^\times$ to
$\C^\times$, $\xi$ can be factored as $\omega \circ \rmn$).
Recall that $\pi = \pi(\chi)$.
Denote uniformizers of $\calo_F$ and $\calo_E$ by
$\varpi_F$ and $\varpi_E$ respectively.

\begin{prop}
	\label{prop:opt-ur}
	Assume that $F / \Q_p$ is unramified if $p = 2$.
	If $\xi$ is unramified,
	then
	\[
		Z\paren{W^\opt, s} =
		\begin{cases}\paren{1 + q^{-s}}
			L\big(\Ad(\pi), s\big)
				&\text{if $E/F$ is not ramified,} \\
			\frac{q^{s}\paren{1 + q^{-s}}}
				{4\paren{1 - q^{-s}}}
				&\text{if $E/F$ is ramified and $p \neq 2$,} \\
			\frac{q^{2s}\paren{1 + q^{-s}}}
				{4\paren{1 - q^{-s}}}
				&\text{if $E/F$ is ramified and $p = 2$.}
		\end{cases}
	\]
	In particular,
	\[
		\calp_\rs\paren{W^\opt}
			= Z\paren{W^\opt, 1}
			=
				\begin{cases}
					\paren{1 + q^{-1}} L\big(\Ad(\pi), 1\big)
						&\text{if $E/F$ is not ramified,} \\
					\frac{q\paren{1 + q^{-1}}}
						{4\paren{1 - q^{-1}}}
						&\text{if $E/F$ is ramified and $p \neq 2$,} \\
					\frac{q^{2}\paren{1 + q^{-1}}}
						{4\paren{1 - q^{-1}}}
						&\text{if $E/F$ is ramified and $p = 2$.}
				\end{cases}
	\]
\end{prop}

\begin{proof}
	First, consider the case that $E/F$ is inert.
	Then $\xi = 1$ and $\delta$ is invertible.
	Let $\eta$ be the quadratic character associated
	to $E/F$.
	The integral reduces to 
	\[
		Z\paren{\Phi^\opt, s}
			= \int_{\verts{t_1} \leq 1} \verts{t_1}^s dt_1
			= \zeta_E(s)
			= \zeta(s) L(\eta, s).
	\]
	Then $L(\Ad(\pi), s) = L(\eta, s)^2 \zeta(s)$,
	and
	\[
		\Psi\paren{\Phi^\opt, s} = 1,
	\]
	(recall the definition of the normalized zeta integral in
	Equation \ref{eq:normalized-zeta-rs}).

	Second, we assume that $E/F$ is split: $E = F \oplus F$.
	Use coordinates $(u_1, u_2)$ for $t_1$ and $(v, v^{-1})$ for $t_2$.
	After composition with a character $\omega \circ \rmn$,
	we may assume that $\chi = (\mu, 1)$ with $\mu$ unramified.
	In this case, $\delta$ is still invertible.
	Therefore
	\begin{align*}
		Z\paren{\Phi^\opt, s}
			&= \int_{\substack{\verts{u_1} \leq 1 \\ \verts{u_2} \leq 1}}
				\verts{u_1 u_2}^s
				\int_{\verts{u_2} \leq \verts{v} \leq \verts{u_1}^{-1}}
				\mu(v) dv du_1 du_2 \\
			&= \sum_{m, n \geq 0} q^{-s(m + n)}
				\sum_{-m \leq \ell \leq n} \mu\paren{\varpi_F^\ell} \\
			&= \frac{1}{1 - \mu(\varpi_F)} \sum_{m, n \geq 0}
				q^{-s(m + n)} \paren{\mu(\varpi_F)^{-m} -\mu(\varpi_F)^{n + 1}}\\
			&= \frac{1}{\big(1 - \mu(\varpi_F)\big) \paren{1 - q^{-s}}}
				\paren{\frac{1}{1 - \mu(\varpi_F)^{-1} q^{-s}} -
				\frac{\mu(\varpi_F)}{1 - \mu(\varpi_F) q^{-s}}} \\
			&= \frac{1 + q^{-s}}
				{\paren{1 - q^{-s}}
				\paren{1 - \mu(\varpi_F)q^{-s}}
				\paren{1 - \mu(\varpi_F)^{-1}q^{-s}}} \\
			&= \paren{1 + q^{-s}}\zeta(s)
				L(\mu, s)L(\mu^{-1}, s).
	\end{align*}
	With our definition,
	$L\paren{\Ad(\pi),s}
	= \zeta(s)L(\mu,s)L(\mu^{-1},s)$.
	Therefore
	\[
		\Psi\paren{\Phi^\opt, s} = 1.
	\]

	Finally, suppose that $E/F$ is ramified.
	Then $\kappa_{E/F} = 1/2$.
	Again, $\xi = \restr{\chi}{E^1} = 1$, and
	\[
		Z\paren{\Phi^\opt,s}
			=
		\frac{1}{2}
		\sum_{n \geq -\ord(\delta)}
			q^{-ns}
			\vol\paren{
				\paren{1 + \varpi_E^{-n}\calo_E}^1}.
	\]

	We can compute the volume case-by-case.
	First, if $n \geq 0$, then $\vol((1 + \varpi^{-n}_E \calo_E)^1) = 1$.
	Second, if $n = -1$, then $E^1$ is the union of
	$\pm (1 + \varpi_E \calo_E)^1$.
	In this case,
	\[
		\vol\paren{\paren{1 + \varpi_E^{1} \calo_E}^1} =
		\begin{cases}
			\frac{1}{2} &\text{if $p \neq 2$,} \\
			1 &\text{if $p = 2$}.
		\end{cases}
	\]
	Finally, we treat the case $n < -1$.
	If $p = 2$ (and $F$ is unramified over $\Q_2$ by assumption),
	then $\ord(\delta) = 2$.
	In this case, $E^1 \pmod{\varpi_E}$ is generated by $1$ and
	$1 + \varpi_E$.
	Thus $\vol((1 + \varpi_E^2 \calo_E)^1) = 1/2$.

	Collecting the terms for all $n$,
	we have the following formula for $Z(\Phi^\opt, s)$
	when $E/F$ is ramified. If $p \neq 2$, then
	\begin{align*}
		Z\paren{\Phi^\opt,s}
			&= \frac{1}{2}
				\paren{
					\frac{1}{2}q^s
					+ \sum_{n \geq 0}q^{-ns}} \\
			&= \frac{q^s\paren{1 + q^{-s}}}
				{4\paren{1 - q^{-s}}}.
	\end{align*}
	If $p = 2$, then
	\begin{align*}
		Z\paren{\Phi^\opt,s}
			&= \frac{1}{2}
				\paren{
					\frac{1}{2}q^{2s}
					+ q^s
					+ \sum_{n \geq 0}q^{-ns}} \\
			&= \frac{q^{2s}\paren{1 + q^{-s}}}
				{4\paren{1 - q^{-s}}}.
	\end{align*}
\end{proof}

\subsection{Ramified calculation: \texorpdfstring{$E/F$}{E/F} inert}

In this subsection, we calculate $\calp_\rs(W^\opt)$ when $\xi$
is ramified and $E/F$ is inert.
\begin{prop}
	\label{prop:opt-inert}
	If $\xi$ is ramified
	and $E/F$ is inert with $p \neq 2$,
	then
	\[
		Z\paren{\Phi^\opt, s}
			= \chi(-1) \paren{1 + q^{-1}}^{-2} q^{2o(\xi)(s - 1)}.
	\]
	In particular,
	\[
		\calp_\rs\paren{W^\opt}
			= Z\paren{\Phi^\opt, 1}
			= \chi(-1) \paren{1 + q^{-1}}^{-2}.
	\]
\end{prop}

\begin{proof} 
	Recall that
	\[
		Z\paren{\Phi^\opt, s}
			= \chi(-1) \int_{\delta^{-1} \calo_c} \verts{t_1}^s 
				\int_{\paren{1 + t_1^{-1}\calo_c}^1} \xi\paren{t_2} dt_2 dt_1.
	\]
	Here $c = c(\xi) = \varpi^{o(\xi)}$ is the
	conductor of $\xi$, with associated order
	$\calo_c = \calo_{o(\xi)}$.
	The element $\varpi$ is a uniformizer of $F$, and
	for every integer $k$, set
	\[
		\calo_k := \calo_F
			+ \varpi^k\calo_F\epsilon.
	\]
	We can write
	$\calo_E = \calo_F + \calo_F\epsilon$, where
	$\epsilon^2 \in \calo_F^\times$.
	Since $p \neq 2$, multiplication by $2$ does not
	change the ideal generated by
	$\varpi^{o(\xi)}\epsilon$. Thus
	\[
		\delta^{-1}\calo_c = \calo_{-o(\xi)}.
	\]
	
	Consider the contribution from
	$t_1 \in \calo_E$.
	The double integral in the expression for
	$Z\paren{\Phi^\opt,s}$ becomes
	\begin{align*}
		\sum_{n \geq 0} q^{-2ns} \int_{\calo_E^\times}
				\int_{\paren{1 + \varpi^{-n} u\calo_c}^1} \xi\paren{t_2} dt_2 du
			&= \sum_{n \geq 0} q^{-2ns} \vol\paren{\calo_c^\times}
				\sum_{u \in \calo_E^\times / \calo_c^\times}
				\int_{\paren{1 + \varpi^{-n} u \calo_c}^1}
				\xi\paren{t_2} dt_2 \\
			&= \vol\paren{\calo_c^\times} \sum_{n \geq 0} q^{-2ns}
				\int_{\paren{1 + \varpi^{-n} \calo_E^\times \calo_c}^1}
				\xi\paren{t_2} dt_2
	\end{align*}
	Notice that
	$\calo_E^\times\calo_c = \calo_E$.
	Thus
	$1 + \varpi^{-n}\calo_E^\times\calo_c
				= 1 + \varpi^{-n}\calo_E
				\supset E^1$,
	so the last integral vanishes. Therefore there is no
	contribution to $Z(\Phi^\opt, s)$ from
	$t_1 \in \calo_E$.

	Now consider the remaining contribution to
	$Z(W^\opt, s)$ from $t_1 \notin \calo_E$.
	For $1 \leq k \leq o(\xi)$, write
	$t_1 = \varpi^{-k}\epsilon u$ with
	$u \in \calo_k^\times$.
	The remaining double integrals become
	\begin{align}
		\label{eq:zeta-inert-notin-COE}
		Z\paren{\Phi^\opt, s}
			&= \chi(-1) \sum_{k = 1}^{o(\xi)} q^{2ks}
				\int_{u \in \calo_k^\times} du
				\int_{\paren{1 + \varpi^k \epsilon u \calo_c}^1}
				\xi\paren{t_2} dt_2 \\
			&= \chi(-1) \vol\paren{\calo_c^\times}
				\sum_{k = 1}^{o(\xi)} q^{2ks}
				\int_{\paren{1 + \varpi^k \epsilon \calo_k^\times \calo_c}^1}
				\xi\paren{t_2} dt_2. \nonumber
	\end{align}

	To further compute the remaining integral, we use the decomposition
	\[
		\calo_c
			= \calo_F + \varpi^{o(\xi)} \calo_E
			= \varpi^{o(\xi)} \calo_E \cup
				\bigcup_{i = 1}^{o(\xi)} \varpi^{o(\xi) - i} \calo_i^\times.
	\]
	Then
	\begin{align*}
		\calo_k^\times \calo_c
			&= \varpi^{o(\xi)} \calo_E \cup
				\bigcup_{i = 1}^k \varpi^{o(\xi) - i} \calo_i^\times \cup
				\bigcup_{j = k + 1}^{o(\xi)} \varpi^{o(\xi) - j} \calo_k^\times \\
			&= \varpi^{o(\xi) - k} \calo_k \cup
				\bigcup_{i = 0}^{o(\xi) - k - 1}
				\paren{\varpi^i \calo_k \setminus \varpi^{i + 1} \calo_{k - 1}}.
	\end{align*}
	For a set $X$, let $\one_X$ denote its characteristic function.
	We have,
	\begin{align*}
		\one_{\paren{1 + \varpi^k \epsilon \calo_k^\times \calo_c}^1}
			&= \one_{\paren{1 + \varpi^{o(\xi)} \epsilon \calo_k}^1}
				+ \sum_{i = 0}^{o(\xi) - k - 1}
				\one_{\paren{1 + \varpi^{k + i} \epsilon \calo_k}^1}
				+ \sum_{j = 0}^{o(\xi) - k - 1}
				\one_{\paren{1 + \varpi^{k + j + 1} \epsilon \calo_{k - 1}}^1} \\
			&= \sum_{i = k}^{o(\xi)}
				\one_{\paren{1 + \varpi^{i} \epsilon \calo_k}^1} +
				\sum_{j = k + 1}^{o(\xi)}
				\one_{\paren{1 + \varpi^{j} \epsilon \calo_{k - 1}}^1}
	\end{align*}
	In particular, the integral becomes
	\begin{align*}
		\int_{\paren{1 + \varpi^k \epsilon \calo_k^\times \calo_c}^1}
				\xi\paren{t_2} dt_2
			&= \sum_{i = k}^{o(\xi)}
				\int_{\paren{1 + \varpi^{i} \epsilon \calo_k}^1}
				\xi\paren{t_2} dt_2 -
				\sum_{j = k + 1}^{o(\xi)}
				\int_{\paren{1 + \varpi^{j} \epsilon \calo_{k - 1}}^1}
				\xi\paren{t_2} dt_2.
	\end{align*}

	Next, we apply the following lemma.
	\begin{lem}
		\label{lem:order-epsilon}
		For any integer $i \geq 0$,
		\[
			\paren{1 + \varpi^i \calo_E}^1 = \paren{1 + \varpi^i \epsilon \calo_i}^1.
		\]
	\end{lem}
	
	\begin{proof}[Proof of Lemma~\ref{lem:order-epsilon}]
		This is trivial when $i = 0$.
		Assume that $i > 0$ and let $x \in (1 + \varpi^i \calo_E)^1$.
		Then we can write for $\alpha, \beta \in \calo_F$,
		\[
			x = 1 + \varpi^{i}\alpha + \varpi^i \beta \epsilon.
		\]
		Factor out $1 + \varpi^i \alpha$ to obtain
		for $\gamma \in \calo_F$,
		\[
			x = \paren{1 + \varpi^i \alpha}
				\frac{1 + \varpi^{i}\alpha + \varpi^i \beta \epsilon}
				{1 + \varpi^i \alpha}
				= \paren{1 + \varpi^i \alpha} \paren{1 + \varpi^i \gamma \epsilon},
		\]
		for some $\gamma \in \calo_F$.
		Now we take the norm $\rmn_{E/F}$ of both sides,
		noting that $1 + \varpi^i \alpha \in F$ and the conjugate of
		$1 + \varpi^i \gamma \epsilon$ is $1 - \varpi^i \gamma \epsilon$,
		\[
			1 = \paren{1 + 2\varpi^{i} \alpha + \varpi^{2i} \alpha^2}
				\paren{1 - \varpi^{2i} \gamma^2 \epsilon^2}.
		\]
		Since $2$ is invertible in $\calo_F$, $\alpha \in \varpi^{i} \calo_F$.
		This shows that $x \in 1 + \varpi^{i} \epsilon \calo_{i}$.
	\end{proof}
	By Lemma~\ref{lem:order-epsilon}, we may replace
		$(1 + \varpi^i \epsilon \calo_k)^1$ and
		$(1 + \varpi^j \epsilon \calo_{k - 1})^1$
		in the integrals by
		$(1 + \varpi^i \calo_E)^1$ and
		$(1 + \varpi^j \calo_E)^1$.
		For $k < o(\xi)$, the nonzero terms in the
		two sums cancel. When $k = o(\xi)$, only the
		first sum contributes, and the integral is
		$\vol\paren{
			1 + \varpi^{o(\xi)}\calo_E}^1$.
		Consequently,
		\[
			Z\paren{\Phi^\opt,s}
				= \chi(-1)\vol\paren{\calo_c^\times}
					q^{2o(\xi)s}
					\vol\paren{
						\paren{
							1 + \varpi^{o(\xi)}
								\calo_E}^1}.
		\]
	To compute the volume of $\calo_E^\times / \calo_c^\times$, we observe
	that $\calo_E^\times$ and $\calo_c^\times$
	both contain $1 + \varpi^{o(\xi)} \calo_E$,
	factor both the top and bottom,
	and then find the cardinality,
	\begin{align*}
		\vol\paren{\calo_E^\times / \calo_c^\times}
			&= \vol\paren{\paren{\calo_E/\varpi^{o(\xi)}}^\times /
				\paren{\calo_F / \varpi^{o(\xi)}}^\times} \\
			&= \frac{\paren{q^2 - 1} q^{2o(\xi) - 2}}
				{(q - 1) q^{o(\xi) - 1}} \\
			&= \paren{1 + q^{-1}} q^{o(\xi)}.
	\end{align*}
	We find the volume of $\calo_E^1 / \paren{1 + \varpi^{o(\xi)} \calo_E}^1$
	in a similar manner.
	\begin{align*}
		\vol\paren{\calo_E^1 / \paren{1 + \varpi^{o(\xi)} \calo_E}^1}
			&= \vol\paren{\paren{\calo_E / \varpi^{o(\xi)}}^1} \\
			&= \vol\paren{\paren{\calo_E / \varpi^{o(\xi)}}^\times /
				\paren{\calo_F / \varpi^{o(\xi)}}^\times} \\
			&= \paren{1 + q^{-1}} q^{o(\xi)}.
	\end{align*}
	Since $E/F$ is unramified, the compatible measures
	satisfy $\vol(\calo_E^\times) = \vol(\calo_E^1) = 1$.
	Therefore
	\[
		Z(\Phi, s) = \chi(-1) \paren{1 + q^{-1}}^{-2} q^{2o(\xi)(s - 1)}.
	\]
\end{proof}

\subsection{Ramified calculation: \texorpdfstring{$E/F$}{E/F} split}
In this subsection, we calculate $\calp_\rs(W^\opt)$ when $\xi$
is ramified and $E/F$ is split. 

\begin{prop}
	\label{prop:opt-split}
	If $\xi$ is ramified
	and $E/F$ is split with $p \neq 2$,
	then
	\[
		Z\paren{\Phi^\opt, s}
			= \chi(-1)\paren{1 - q^{-1}}^{-2}
				\paren{\frac{2q^{-2o(\xi) - s}}{1 - q^{-s}}
				+ q^{2o(\xi)(s - 1)}}.
	\]
	In particular,
	\[
		\calp_\rs\paren{W^\opt}
			= Z\paren{W^\opt, 1}
			= \chi(-1)\paren{1 - q^{-1}}^{-3}
				\paren{1 - q^{-1} + 2q^{-2o(\xi) - 1}}.
	\]
\end{prop}

\begin{proof}
	Recall that
	\[
		Z\paren{\Phi^\opt, s}
			= \chi(-1)
				\int_{\delta^{-1}\calo_c}
					\verts{t_1}^s
				\int_{\paren{
					1 + t_1^{-1}\calo_c}^1}
					\xi\paren{t_2}dt_2dt_1.
	\]
	We have
	$\calo_c
	= \calo_F + \varpi^{o(\xi)}\calo_F\epsilon$
	and
	$\delta = \varpi^{o(\xi)}\epsilon$.
	Therefore
	\[
		\delta^{-1}\calo_c
			= \calo_{-o(\xi)}
			= \calo_F
				+ \varpi^{-o(\xi)}
					\epsilon\calo_F.
	\]
	Let us first consider the contribution to
	$Z(\Phi^\opt, s)$
	from the integral over $t_1 \in \calo_E$,
	\[
		\chi(-1) \vol\paren{\calo_c^\times}
			\sum_{m, n \geq 0} q^{-(m + n)s}
			\int_{\paren{1 + \paren{\varpi^{-m}, \varpi^{-n}} \calo_E^\times \calo_c}^1}
			\xi\paren{t_2} dt_2.
	\]
	
	We can decompose $\calo_c$ as
	\[
		\calo_c
			= \calo_F + \varpi^{o(\xi)} \calo_E
			= \varpi^{o(\xi)} \calo_E \cup \bigcup_{k = 0}^{o(\xi) - 1}
				\varpi^k \calo_{o(\xi) - k}^\times.
	\]
	For a set $X$, let $\one_X$ denote its characteristic function.
	Then
	\begin{align*}
	\one_{\calo_E^\times \calo_c}
		&= \one_{\varpi^{o(\xi)} \calo_E}
			+ \sum_{k = 0}^{o(\xi) - 1} \one_{\varpi^k \calo_E^\times} \\
		&= \one_{\varpi^{o(\xi)} \calo_E}
			+ \sum_{k = 0}^{o(\xi) - 1}
			\paren{\one_{\varpi^k\calo_E}
				- \one_{\paren{\varpi^{k + 1},\varpi^k}\calo_E}
				- \one_{\paren{\varpi^k,\varpi^{k + 1}}\calo_E}
				+ \one_{\paren{\varpi^{k + 1},\varpi^{k + 1}}\calo_E}} \\
		&= \one_{\calo_E}
			+ 2\sum_{i = 1}^{o(\xi)} \one_{\varpi^i \calo_E}
			- \sum_{j = 0}^{o(\xi) - 1}
			\paren{\one_{\paren{\varpi^{j + 1}, \varpi^j} \calo_E}
			+ \one_{\paren{\varpi^j, \varpi^{j + 1}} \calo_E}}.
	\end{align*}
	In particular, the integral becomes
	\begin{align}
		\label{eq:integral-decomp-ramified-split}
		\int_{\paren{1 + \paren{\varpi^{-m}, \varpi^{-n}} \calo_E^\times \calo_c}^1}
			\xi\paren{t_2} dt_2
			= \, &\int_{\paren{\varpi^{-m}, \varpi^{-n}} \calo_E}
				\xi\paren{t_2} dt_2
				+ 2\sum_{i = 1}^{o(\xi)}
				\int_{\paren{\varpi^{i - m}, \varpi^{i - n}} \calo_E}
				\xi\paren{t_2} dt_2 \\
				\, &- \sum_{j = 0}^{o(\xi) - 1}
				\paren{\int_{\paren{\varpi^{j + 1 - m}, \varpi^{j - n}} \calo_E}
				\xi\paren{t_2} dt_2
				+ \int_{\paren{\varpi^{j - m}, \varpi^{j + 1 - n}} \calo_E}
				\xi\paren{t_2} dt_2}. \nonumber
	\end{align}
	So we need to compute, for integers $a$ and $b$,
	\[
		\int_{\paren{1 + \paren{\varpi^a, \varpi^b} \calo_E}^1} \xi(t) dt.
	\]
	
	\begin{lem}
		\label{lem:integral-ab}
		\[
			\int_{\paren{1 + \paren{\varpi^a, \varpi^b} \calo_E}^1} \xi(t) dt = 
				\begin{cases}
					q^{-\max(a, b)} \paren{1 - q^{-1}}^{-1}
						&\text{if $\max(a, b) \geq o(\xi)$}, \\
					0 &\text{otherwise}.
				\end{cases}
		\]
	\end{lem}
	
	\begin{proof}[Proof of Lemma~\ref{lem:integral-ab}]
		Using the coordinates $(t, t^{-1})$, the integral is then over
		$t$ such that
		\begin{align*}
			t &= 1 + \varpi^a x, \\
			t^{-1} &= 1 + \varpi^b y,
		\end{align*}
		for $x, y \in \calo_F$.
		
		If $\max(a,b) \leq 0$, then these conditions become
		\[
			\verts{\varpi}^{-b}
				\leq \verts{t}
				\leq \verts{\varpi}^{a}.
		\]
		Hence the set of such $t$ is stable under
		multiplication by $\calo_F^\times$.
		The integral vanishes because $\xi$ is ramified.

		If $\max(a,b) > 0$, then the integral is
		\[
			\int_{\paren{
				1 + \varpi^{\max(a,b)} \calo_E}^1}
				\xi(t)dt
				=
			\begin{cases}
				\vol\paren{
					1 + \varpi^{\max(a, b)} \calo_F}
					&\text{if $\max(a, b) \geq o(\xi)$}, \\
				0
					&\text{if
						$0 < \max(a, b) < o(\xi)$}.
			\end{cases}
		\]
		Evaluating the volume proves the lemma.
	\end{proof}

	By Lemma~\ref{lem:integral-ab}, the terms of
	Equation~\ref{eq:integral-decomp-ramified-split} either individually
	vanish or cancel completely with each other,
	unless either $m = 0$ and $n > 0$ or $m > 0$ and $n = 0$.
	In particular,
	\[
		\int_{\paren{1 + \paren{\varpi^{-m}, \varpi^{-n}} \calo_E^\times \calo_c}^1}
			\xi\paren{t_2} dt_2 =
			\begin{cases}
				q^{-o(\xi)} \paren{1 - q^{-1}}^{-1}
					&\text{if $mn = 0$ and $m + n > 0$}, \\
				0 &\text{otherwise}.
			\end{cases}
	\]
	Then the sum over all such $m$ and $n$ is twice
	the geometric series in $q^{-s}$ from fixing $m$ or $n$ to be zero,
	\[
		\sum_{m, n \geq 0} q^{-(m + n)s}
		\int_{\paren{1 + \paren{\varpi^{-m}, \varpi^{-n}} \calo_E^\times \calo_c}^1}
		\xi\paren{t_2} dt_2
		= 2 \frac{q^{-s}}{1 - q^{-s}} q^{-o(\xi)} \paren{1 - q^{-1}}^{-1}.
	\]
	The remaining calculation in the $t_1 \in \calo_E$ case is
	$\vol(\calo_c^\times)$, which was done in the inert case
	(cf. proof of Proposition~\ref{prop:opt-inert}),
	\begin{align*}
		\vol\paren{\calo_E^\times / \calo_c^\times}
			&= \vol\paren{\paren{\calo_E/\varpi^{o(\xi)}}^\times /
				\paren{\calo_F / \varpi^{o(\xi)}}^\times} \\
			&= \frac{\paren{q^2 - 1} q^{2o(\xi) - 2}}
				{(q - 1) q^{o(\xi) - 1}} \\
			&= \paren{1 + q^{-1}} q^{o(\xi)}.
	\end{align*}
	Therefore the contribution from $t_1 \in \calo_E$ to
	$Z(W^\opt, s)$ is,
	\begin{align*}
		\chi(-1) \vol\paren{\calo_{o(\xi)}^\times} \sum_{m, n \geq 0}
			q^{-(m + n)s}
			\int_{\paren{1 + \paren{\varpi^{-m}, \varpi^{-n}} \calo_E^\times \calo_c}^1}
			&\xi\paren{t_2} dt_2 \\
			&= 2\chi(-1) \paren{1 - q^{-1}}^{-2}
			\frac{q^{-2o(\xi) - s}}{1 - q^{-s}}.
	\end{align*}

	Now we consider the remaining contribution to $Z(W^\opt, s)$
	from $t_1 \notin \calo_E$.
	We use,
	\[
		\calo_c = \varpi^{o(\xi)} \calo_E \cup
			\bigcup_{k = 1}^{o(\xi)} \varpi^{o(\xi) - k} \calo_k^\times,
	\]
	to get,
	\[
		\delta^{-1} \calo_c \setminus \calo_E
			= \bigcup_{k = 1}^{o(\xi)} \varpi^{-k} \epsilon \calo_{k}^\times.
	\]
	Then we have that the $t_1 \notin \calo_E$
	contribution to $Z(W^\opt, s)$ is,
	\begin{align*}
		\chi(-1) \int_{\delta^{-1} \calo_c \setminus \calo_E}
			\verts{t_1}^s 
			\int_{\paren{1 + t_1^{-1} \calo_c}^1} \xi\paren{t_2} dt_2 dt_1
			&= \chi(-1) \sum_{k = 1}^{o(\xi)} q^{2ks}
			\vol\paren{\calo_c^\times}
			\int_{\paren{1 + \varpi^k \epsilon \calo_k^\times \calo_c}^1}
			\xi\paren{t_2} dt_2.
	\end{align*}
	This is the same expression from
	Equation~\ref{eq:zeta-inert-notin-COE} that we calculated
	in the proof of Proposition~\ref{prop:opt-inert},
	\begin{align*}
		 \chi(-1) \sum_{k = 1}^{o(\xi)} q^{2ks}
			\vol\paren{\calo_c^\times}
			\int_{\paren{1 + \varpi^k \epsilon \calo_k^\times \calo_c}^1}
			\xi\paren{t_2} dt_2
			= \chi(-1) \paren{1 - q^{-1}}^{-2} q^{2o(\xi)(s - 1)}
	\end{align*}

	Combining the contributions from $t_1 \in \calo_E$ and
	$t_1 \in \delta^{-1} \calo_c \setminus \calo_E$,
	we have shown that
	\[
		Z(\Phi^\opt, s) = \chi(-1) \paren{1 - q^{-1}}^{-2}
			\paren{\frac{2q^{-2o(\xi) - s}}{1 - q^{-s}} +
			q^{2o(\xi)(s-1)}}.
	\]
\end{proof}


\section{Local Rankin--Selberg periods of newforms}
\label{sec:inner-product-new}

The numerator of the local factor in
Theorem~\ref{thm:ratio} is the Rankin--Selberg period
of the new vectors. Let $\Phi_1$ and $\Phi_2$ be the
standard Schwartz functions for $\chi$ and $\chi^{-1}$
from Section~\ref{sec:quadratic-cases}, and set
$\Phi := \Phi_1 \otimes \Phi_2$.
The associated theta-Whittaker functions give the new
vectors when $E/F$ is not ramified. In the ramified
field case, the new vectors are sums of two theta
functions, as described before
Proposition~\ref{prop:theta-new}.

We again view the antinorm $\xi = \chi^{1 - c_{E/F}}$ as the
restriction of $\chi^{-1}$ to $E^1$. By
Proposition~\ref{prop:rs-phi}, the calculation reduces to
\[
	Z\paren{\Phi, s}
		= \int_{E^\times} \verts{t_1}^s
			\int_{E^1}
				\wt{\Phi}\paren{t_1\paren{1 + t_2j}}
				\xi(t_2) dt_2 dt_1,
\]
where
\[
	\wt{\Phi}
		= \int_{\SL_2(\calo_F)} r(k)\Phi\,dk.
\]
In each ramification case, we compute the Fourier
transform $\wh{\Phi}$, determine a congruence subgroup
fixing $\Phi$, use a coset decomposition to calculate
$\wt{\Phi}$, and then evaluate this integral.

\subsection{Unramified calculation}
First, we assume that $\chi$ is unramified.
We do not treat the $E/F$ ramified case
when $p = 2$.
\begin{prop}
	\label{prop:new-ur}
	If $\chi$ is unramified, then
	$\calp_\rs(W^\new) \in \Q(\xi + \xi^{-1})^\times$.
	Furthermore,
	\[
		Z\paren{W^\new, s} =
			\begin{cases}
				Z\paren{W^\opt, s}
					&\text{if $E/F$ is not ramified}, \\
				4(q + 1)^{-1}Z\paren{W^\opt, s}
					&\text{if $E/F$ is ramified and } p \neq 2.
			\end{cases}
	\]
	In particular,
	\[
		\calp_\rs\paren{W^\new} =
			\begin{cases}
				\calp_\rs\paren{W^\opt}
					&\text{if $E/F$ is not ramified}, \\
				4(q + 1)^{-1} \calp_\rs\paren{W^\opt}
					&\text{if $E/F$ is ramified and } p \neq 2.
			\end{cases}
	\]
\end{prop}

\begin{proof}
	With $\Phi_i$ as the characteristic function of $\calo_E$
	(cf. Section~\ref{sec:quadratic-cases}),
	$\Phi$ is the characteristic function of
	$\calo_B^\new := \calo_E + \calo_E j$ in the quaternion
	algebra $B = E + Ej$.
	
	If $E/F$ is unramified, then $\Phi = \Phi^\opt$
	and the result follows from Proposition \ref{prop:opt-ur}.

	Assume that $E/F$ is ramified. We first compute
	$\wt{\Phi}$. Let $d_E\calo_F$ be the discriminant
	ideal of $E/F$. The function $\Phi$ is invariant under
	$B(\calo_F)$. We claim that it is invariant under
	$U_0(d_E)$.

	Let $\delta_{E/F}$ be the different of $E/F$.
	With the self-dual measure used in the Weil
	representation, the dual lattice of $\calo_E$ is
	$\delta_{E/F}^{-1}\calo_E$ and
	$\vol(\calo_E)
	= \verts{\delta_{E/F}}_E^{\frac{1}{2}}$.
	Therefore
	\[
		\wh{\Phi}_i
			= \verts{\delta_{E/F}}_E^{\frac{1}{2}}
				\restr{\one}{
					\delta_{E/F}^{-1}\calo_E},
	\]
	and consequently
	\[
		\wh{\Phi}
			= \verts{d_E}_F
				\restr{\one}{
					\delta_{E/F}^{-1}
					\calo_B^\new}.
	\]
	It follows that $\wh{\Phi}$ is invariant under
	$N(d_E\calo_F)$. Since $U_0(d_E)$ is generated by
	$B(\calo_F)$ and $wN(d_E\calo_F)w^{-1}$, the claim
	follows.

	Since $p \neq 2$, the discriminant ideal is the maximal
	ideal of $\calo_F$; write
	$d_E\calo_F = \varpi_F\calo_F$. The right-coset
	decomposition
	\[
		\SL_2(\calo_F)
			= U_0(\varpi_F)
				\sqcup
				\bigsqcup_{
					b \in \calo_F/\varpi_F\calo_F}
					n(b)wU_0(\varpi_F)
	\]
	and $\vol(\SL_2(\calo_F)) = 1$ give
	\[
		\wt{\Phi}
			= \frac{1}{q + 1}
				\paren{
					\Phi
					+ \sum_{
						b \in
						\calo_F/\varpi_F\calo_F}
						n(b)\wh{\Phi}}.
	\]

	For $x, y \in \calo_E$, the character sum in the
	second term is nonzero precisely when
	$\rmn_{E/F}(x) - \rmn_{E/F}(y)$ belongs to
	$\varpi_F\calo_F$. Since $E/F$ is ramified and
	$p \neq 2$, this is equivalent to
	$x \equiv \pm y \pmod{\varpi_E}$. Set
	\[
		\calo_B^\pm
			:= \calo_E + \calo_Ej
				+ \frac{1 \pm j}{\varpi_E}\calo_E,
	\]
	and let $\Phi^\pm$ be their characteristic functions.
	One has
	$\calo_B^- = \End_{\calo_F}(\calo_E)$ and
	$\calo_B^+
	= \varpi_E\calo_B^-\varpi_E^{-1}$.
	Thus both orders are maximal and contain $\calo_E$
	optimally, while their intersection is
	$\calo_B^\new$. Therefore
	\[
		\sum_{b \in \calo_F/\varpi_F\calo_F}
			n(b)\wh{\Phi}
				= \Phi^+ + \Phi^- - \Phi,
	\]
	and hence
	\[
		\wt{\Phi}
			= \frac{1}{q + 1}
				\paren{\Phi^+ + \Phi^-}.
	\]

	The calculation in
	Proposition~\ref{prop:opt-ur} applies with either
	sign, equivalently after replacing $j$ by $-j$.
	Linearity therefore gives
	\[
		Z\paren{\Phi, s}
			= \frac{2}{q + 1}
				Z\paren{W^\opt, s}.
	\]

	The two summands in each ramified newvector are
	supported on the two cosets of $\GL_2(F)^+$ in
	$\GL_2(F)$. After expanding the zeta integral, the
	two terms with incompatible supports vanish. The two
	remaining terms are equal after right translation by
	$a(\epsilon)$. Since $\verts{\epsilon} = 1$, this
	translation preserves both $f(g, s)$ and the quotient
	measure. Thus
	\[
		Z\paren{W^\new, s}
			= 2Z\paren{\Phi, s}
			= \frac{4}{q + 1}
				Z\paren{W^\opt, s}.
	\]
	Setting $s = 1$ gives the asserted period identity.
	Its rationality and non-vanishing follow from
	Corollary~\ref{cor:opt-rat}.
\end{proof}

\subsection{Ramified calculation: \texorpdfstring{$E/F$}{E/F} inert}
Next, we assume that $\chi$ is ramified and $E/F$ is inert.

\begin{prop}
	\label{prop:new-inert}
	If $\chi$ is ramified and $E/F$ is inert, then
	\[
		\calp_\rs\paren{W^\new}
			= Z\paren{\Phi, 1}
			=	\begin{cases}
				1
					&\text{if $\xi^2$ is unramified} \\
				\frac{1}{q + 1}
					&\text{if $\xi^2$ is ramified}.
			\end{cases}
	\]
\end{prop}

\begin{proof}
	Let $o(\chi)$ be the conductor exponent of $\chi$,
	namely the least positive integer such that $\chi$ is
	trivial on $1 + \varpi_E^{o(\chi)}\calo_E$.
	Since $\chi$ is ramified, it is nontrivial on
	$1 + \varpi_E^{o(\chi)-1}\calo_E$.
	Set $c(\chi) := \varpi_E^{o(\chi)}$.
	Again, the strategy to evaluate $Z(\Phi, 1)$ is to use
	Proposition~\ref{prop:rs-phi} and
	a description of $\wt{\Phi}$ (also using
	a description of $\wh{\Phi}$).

	Since $E/F$ is inert,
	the $\Phi_i$ are the
	restrictions of $\chi^{-1}$ and $\chi$
	on $\calo_E^\times$ respectively
	(cf. Section~\ref{sec:quadratic-cases}).
	Then $\Phi$ is invariant under $B(\calo_F)$.
	Thus it is invariant under some $U_0(\varpi^k)$ for some $k$, which we
	call the level of $\Phi$. To determine such $k$, let us compute
	$\wh{\Phi} = \wh{\Phi}_1 \otimes \wh{\Phi}_2$:
	\begin{align*}
		\wh{\Phi}_1(x) &= \int_{\calo_E^\times} \chi^{-1}(u) \psi_E(xu) du, \\
		\wh{\Phi}_2(y) &= \int_{\calo_E^\times} \chi(v) \psi_E(-yv) dv,
	\end{align*}
	where the measure is additive so that $\vol(\calo_E) = 1$
	and $\psi_E = \psi \circ \Tr_{E/F}$.
	These are local Gauss integrals. Define
	\begin{equation}
		\label{eq:def-epsilon}
		\epsilon(\chi,\psi)
			:= \int_{\calo_E^\times}
				\chi\paren{\varpi_E^{-o(\chi)}u}
				\psi_E\paren{\varpi_E^{-o(\chi)}u}du.
	\end{equation}

	\begin{lem}
		\label{lem:epsilon}
		Let $\chi: E^\times \lra \C^\times$ be a
		ramified multiplicative character with conductor
		exponent $o(\chi) > 0$, and let
		$\psi_E: E \lra \C^\times$ be an additive character
		with conductor exponent $0$. Then
		\begin{enumerate}
			\item
				\[
					\int_{\calo_E^\times}
						\chi(u)\psi_E(xu)du
						=
						\begin{cases}
							\chi^{-1}(x)
								\epsilon(\chi,\psi)
								&\text{if }
									x \in
									\varpi_E^{-o(\chi)}
									\calo_E^\times, \\
							0
								&\text{otherwise;}
						\end{cases}
				\]
			\item
				\[
					\epsilon(\chi,\psi)
					\epsilon(\chi^{-1},\psi^{-1})
						= \verts{
							\varpi_E^{o(\chi)}}_E.
				\]
		\end{enumerate}
	\end{lem}

	\begin{proof}[Proof of Lemma~\ref{lem:epsilon}]
		Put $r := o(\chi)$, and define
		$H_0 := \calo_E^\times$ and
		$H_i := 1 + \varpi_E^i\calo_E$ for $i \geq 1$.
		If $u \in \calo_E^\times$ and
		$v \in \varpi_E^r\calo_E$, then
		$u + v = u\paren{1 + u^{-1}v}$
		and $\chi(u + v) = \chi(u)$ by the definition of
		$r$. Thus the integral over each additive coset
		of $\varpi_E^r\calo_E$ is a scalar multiple of
		$\int_{\varpi_E^r\calo_E}\psi_E(xv)dv$.
		Additive-character orthogonality shows that this
		integral vanishes unless
		$x \in \varpi_E^{-r}\calo_E$.

		Suppose that
		$x \in \varpi_E^{1-r}\calo_E$.
		For every $h \in H_{r-1}$ and
		$u \in \calo_E^\times$, one has
		$x u(h - 1) \in \calo_E$, so
		$\psi_E(xuh) = \psi_E(xu)$.
		The change of variables $u \mapsto uh$ therefore
		multiplies the integral by $\chi(h)^{-1}$.
		Since $\chi$ is nontrivial on $H_{r-1}$, the
		integral vanishes. Its support is consequently
		$\varpi_E^{-r}\calo_E^\times$.

		If $x = \varpi_E^{-r}a$ with
		$a \in \calo_E^\times$, then the change of
		variables $u \mapsto a^{-1}u$ gives
		\begin{align*}
			\int_{\calo_E^\times}
				\chi(u)\psi_E(xu)du
				&= \chi^{-1}(a)
					\int_{\calo_E^\times}
						\chi(u)
						\psi_E(\varpi_E^{-r}u)du \\
				&= \chi^{-1}(x)
					\epsilon(\chi,\psi).
		\end{align*}
		This proves the first identity.

		For the second identity, changing variables
		$u = vw$ gives
		\begin{align*}
			\epsilon(\chi,\psi)
			\epsilon(\chi^{-1},\psi^{-1})
				= \int_{\calo_E^\times}\chi(w)
					\int_{\calo_E^\times}
						\psi_E\paren{
							\varpi_E^{-r}(w - 1)v}
					dv\,dw.
		\end{align*}
		Since
		$\calo_E^\times
		= \calo_E - \varpi_E\calo_E$,
		additive-character orthogonality gives
		\[
			\int_{\calo_E^\times}
				\psi_E\paren{
					\varpi_E^{-r}(w - 1)v}dv
				= \one_{H_r}(w)
					- \verts{\varpi_E}_E
						\one_{H_{r-1}}(w).
		\]
		It follows that
		\begin{align*}
			\epsilon(\chi,\psi)
			\epsilon(\chi^{-1},\psi^{-1})
				&= \int_{H_r}\chi(w)dw
					- \verts{\varpi_E}_E
						\int_{H_{r-1}}\chi(w)dw \\
				&= \vol(H_r) \\
				&= \verts{\varpi_E^r}_E.
		\end{align*}
		Here $\chi$ is trivial on $H_r$ and nontrivial
		on $H_{r-1}$, while multiplication by a unit
		preserves the additive Haar measure.
	\end{proof}

	By Lemma~\ref{lem:epsilon}, we have the following,
	\begin{align*}
		w\Phi 
			&= \wh{\Phi} \\
			&= \restr{\epsilon(\chi, \psi) \epsilon(\chi^{-1}, \psi^{-1})
				\paren{\chi \otimes \chi^{-1}}}{\varpi^{-o(\chi)}
				\paren{\calo_E^\times + \calo_E^\times j}} \\
			&= \restr{q^{-2o(\chi)} \paren{\chi \otimes \chi^{-1}}}
				{\varpi^{-o(\chi)} \paren{\calo_E^\times +\calo_E^\times j}}.
	\end{align*}
	It follows that $\Phi$ is invariant under $U_0(\varpi^{2o(\chi)})$.
	So we can take $k = 2o(\chi)$.
	
	Now we calculate $\wt{\Phi}$ using the Bruhat decomposition,
	\[
		\SL_2\paren{\calo_F}
			= w N\paren{\calo_F / \varpi^k} w U_0 \paren{\varpi^k} \, \cup \,
			N\paren{\varpi \calo_F / \varpi^k} w U_0 \paren{\varpi^k},
	\]
	so
	\[
		\wt{\Phi}(x)
			= \paren{q^{k} + q^{k - 1}}^{-1}
				\paren{w \sum_{b \in \varpi\calo_F / \varpi^k} r(b) w\Phi(x) +
				\sum_{b \in \calo_F / \varpi^k} r(b) w\Phi(x)}.
	\]
	The two sums are respectively equal to:
	\begin{align*}
		\sum_{b \in \calo_F / \varpi^k} r(b) \wh{\Phi}(x + yj)
			&= \sum_{b \in \calo_F / \varpi^k} \psi\big(b(\rmn(x) - \rmn(y))\big)
				\wh{\Phi}(x + yj). \\
		\sum_{b \in \varpi\calo_F / \varpi^k} r(b) \wh{\Phi}(x + yj)
			&= \sum_{b \in \varpi\calo_F / \varpi^k} \psi\big(b(\rmn(x) - \rmn(y))\big)
				\wh{\Phi}(x + yj).
	\end{align*}
	These sums are non-zero only if $\rmn(x) - \rmn(y) \in \calo_F$ and
	$\rmn(x) - \rmn(y) \in \varpi^{-1} \calo_F$, respectively.
	Write $x = \varpi^{-k/2} u$ and $y = \varpi^{-k/2}v$ with
	$u, v\in \calo_E^\times$. 
	Then these conditions are that
	$\rmn(uv^{-1}) \in 1 + \varpi^k \calo_F$ and
	$\rmn(uv^{-1}) \in 1 + \varpi^{k - 1} \calo_F$,
	respectively. So,
	\[
		\wt{\Phi}
			= q^{-k} (q + 1)^{-1}
				\paren{\restr{w\paren{\chi \otimes \chi^{-1}}}{\Omega_1} +
				q \paren{\restr{\chi \otimes \chi^{{-1}}}{\Omega_0}}},
	\]
	where for integers $i \in \set{0, 1}$,
	\[
		\Omega_{i}
			:= \set{(u, v) \in \varpi^{-\frac{k}{2}}
				\paren{\calo_E^\times \times \calo_E^\times} \Mid
				\rmn\paren{\frac{u}{v}} \in 1 + \varpi^{k - i} \calo_F}.
	\]
	To further describe $\wt{\Phi}$, we expand the
	congruence condition defining $\Omega_1$ by finite
	characters. Set
	\begin{align*}
		A_{k-1}
			&:= \paren{\calo_F/\varpi^{k-1}\calo_F}^\times,
		\\
		\wh{A}_{k-1}
			&:= \Hom\paren{A_{k-1},\C^\times}.
	\end{align*}
	For $\omega \in \wh{A}_{k-1}$, inflate $\omega$
	to $\calo_F^\times$, extend it to $F^\times$ by
	$\omega(\varpi) = 1$, and set
	$\chi_\omega := \chi\cdot\paren{\omega\circ\rmn_{E/F}}$.
	For every integer $n$ and $x \in \calo_F^\times$,
	define
	\[
		f_n(x)
			:= \verts{A_{k-1}}^{-1}
				\sum_{\substack{
					\omega \in \wh{A}_{k-1} \\
					0 < o(\chi_\omega) = k/2 + n}}
				\omega(x).
	\]
	The index $n$ records the conductor exponent of the
	twist $\chi_\omega$, which determines the support of
	its Gauss integral.
	Define $\Phi_n$ by
	\[
		\Phi_n(x,y)
			:= \begin{cases}
				q^{k-2n}\chi\paren{y/x}
					f_n\paren{\rmn(y/x)}
					&\text{if $(x,y) \in
						\varpi^{-n}(\calo_E^\times)^2$,} \\
				0
					&\text{otherwise.}
			\end{cases}
	\]
	Finally, set
	\begin{align*}
		J(x)
			&:= \one_{\calo_E}(x)
				- q^{-2}
					\one_{\varpi_E^{-1}\calo_E}(x), \\
		\Phi^{\mathrm{ur}}(x,y)
			&:= \verts{A_{k-1}}^{-1}q^{2k}
				J\paren{\varpi^{-k/2}x}
				J\paren{\varpi^{-k/2}y}.
	\end{align*}

	\begin{lem}
		\label{lem:newform-wtphi-inert}
		With the above notation,
		\[
			\restr{w\paren{\chi\otimes\chi^{-1}}}
				{\Omega_1}
				= \sum_n\Phi_n
					+ \begin{cases}
						\Phi^{\mathrm{ur}}
							&\text{if $\xi = 1$,} \\
						0
							&\text{otherwise.}
					\end{cases}
		\]
	\end{lem}

	\begin{proof}[Proof of
		Lemma~\ref{lem:newform-wtphi-inert}]
		After the change of variables
		$(u,v) \mapsto
		(\varpi^{-k/2}u,\varpi^{-k/2}v)$,
		the set $\Omega_1$ becomes the subgroup
		\[
			G
				:= \Ker\paren{
					\begin{aligned}
						(\calo_E^\times)^2
							&\lra A_{k-1} \\
						(u,v)
							&\longmapsto
								\rmn\paren{v/u}.
					\end{aligned}}
		\]
		Consequently,
		\begin{align*}
			\restr{w\paren{\chi\otimes\chi^{-1}}}
				{\Omega_1}(x,y)
				&= q^{2k}\int_G\chi\paren{u/v}
					\psi_E\paren{
						\varpi^{-k/2}xu
						- \varpi^{-k/2}yv}
					du\,dv.
		\end{align*}
		Character orthogonality gives
		\[
			\one_G(u,v)
				= \verts{A_{k-1}}^{-1}
					\sum_{\omega\in\wh{A}_{k-1}}
						\omega\paren{\rmn(u/v)}.
		\]
		The contribution of $\omega$ is therefore the
		product
		\begin{align*}
			&\int_{\calo_E^\times}\chi_\omega(u)
				\psi_E\paren{\varpi^{-k/2}xu}du
			\, \cdot \,
			\int_{\calo_E^\times}\chi_\omega^{-1}(v)
				\psi_E\paren{-\varpi^{-k/2}yv}dv.
		\end{align*}

		Suppose first that $\chi_\omega$ is ramified.
		By Lemma~\ref{lem:epsilon}, this product is nonzero
		only if
		\[
			\ord(x) = \ord(y)
				= k/2-o(\chi_\omega).
		\]
		If this common value is $-n$, then the product is
		$q^{-2o(\chi_\omega)} \chi_\omega\paren{y/x}$.
		Multiplying by $q^{2k}$ and summing over the
		characters with $o(\chi_\omega)=k/2+n$ gives
		$\Phi_n(x,y)$.

		Since $E/F$ is unramified, the norm map
		$\calo_E^\times \lra \calo_F^\times$ is
		surjective with kernel $E^1$.
		Thus a twist $\chi_\omega$ is unramified precisely
		when $\xi = 1$, and in that case there is exactly
		one such $\omega$. For an unramified character,
		its restriction to $\calo_E^\times$ is trivial, and
		\[
			\int_{\calo_E^\times}\psi_E(zu)du
				= J(z).
		\]
		Its contribution is therefore
		$\Phi^{\mathrm{ur}}(x,y)$. This proves the lemma.
	\end{proof}

	Set
	\[
		\Phi_{\Omega_1}
			:= \sum_n\Phi_n
				+ \begin{cases}
					\Phi^{\mathrm{ur}}
						&\text{if $\xi = 1$,} \\
					0
						&\text{otherwise.}
				\end{cases}
	\]
	By Lemma~\ref{lem:newform-wtphi-inert},
	\[
		\wt{\Phi}
			= q^{-k}(q + 1)^{-1}
				\paren{
					\Phi_{\Omega_1}
					+ q\restr{\chi\otimes\chi^{-1}}
						{\Omega_0}}.
	\]
	We now evaluate the terms at $s = 1$ using
	Proposition~\ref{prop:rs-phi}. First,
	\begin{align*}
		Z\paren{\Phi_n,1}
			&= q^{k-2n}
				\int_{\varpi^{-n}\calo_E^\times}
					\verts{t_1}dt_1
				\int_{E^1}f_n(1)dt_2 \\
			&= q^k f_n(1),
	\end{align*}
	where we used
	$\restr{\xi}{E^1}
	= \restr{\chi^{-1}}{E^1}$.

	For the unramified term, the two possible values of
	$J(\varpi^{-k/2}t_1)$ give
	\begin{align*}
		Z\paren{\Phi^{\mathrm{ur}},1}
			&= \frac{q^{2k}}{\verts{A_{k-1}}}
				\paren{
					q^{-k + 2}q^{-4}
					+ (1 - q^{-2})^2
						\sum_{m\geq k/2}q^{-2m}} \\
			&= \frac{q^k}{\verts{A_{k-1}}}.
	\end{align*}
	Finally,
	\[
		Z\paren{
			\restr{\chi\otimes\chi^{-1}}{\Omega_0},1}
			= \begin{cases}
				q^k
					&\text{if $\xi^2 = 1$,} \\
				0
					&\text{otherwise.}
			\end{cases}
	\]
	Indeed, the inner integral over $E^1$ is the integral
	of $\xi^2$.

	The characters in $\wh{A}_{k-1}$ are partitioned by
	the conductor exponents of the ramified
	$\chi_\omega$, with one unramified character omitted
	exactly when $\xi = 1$. Hence
	\[
		\sum_n f_n(1)
			= \begin{cases}
				1-\verts{A_{k-1}}^{-1}
					&\text{if $\xi = 1$,} \\
				1
					&\text{otherwise.}
			\end{cases}
	\]
	Combining the three contributions gives
	\[
		Z\paren{\Phi,1}
			= \begin{cases}
				1
					&\text{if $\xi^2 = 1$,} \\
				\frac{1}{q + 1}
					&\text{otherwise.}
			\end{cases}
	\]
	Since $\xi(\varpi_E)=1$, the character $\xi^2$ is
	unramified precisely when $\xi^2 = 1$. This proves
	Proposition \ref{prop:new-inert}.
\end{proof}

\subsection{Fully ramified calculation:
	\texorpdfstring{$E/F$}{E/F} split}
Next, we assume that $\chi$ is ramified and $E/F = F \oplus F$ is
split. In particular, we can write $\chi = (\chi_1, \chi_2)$.
There are two cases: either both $\chi_1$ and $\chi_2$ are ramified,
or exactly one of them is ramified.
Here, we consider the ``fully ramified'' case that both
$\chi_i$ are ramified.

\begin{prop}
	\label{prop:new-fr}
	If $\chi = (\chi_1, \chi_2)$ with $\chi_1$ and $\chi_2$ ramified
	and $E = F \oplus F$, then
	\[
		\calp_\rs\paren{W^\new}
			= Z\paren{\Phi, 1}
			= \begin{cases}
					1 &\text{if $\xi^2$ is unramified}, \\ 
					\frac{1}{q + 1} &\text{if $\xi^2$ is ramified}.
				\end{cases}
	\]
\end{prop}

\begin{proof}
	Let $o_1, o_2$ be the orders of $\chi_1$ and $\chi_2$ respectively.
	Take 
	$\Phi_\chi$ to be $\chi^{-1}$ restricted to $\calo_E^\times$, and
	$\Phi = \Phi_\chi \otimes \Phi_{\chi^{-1}}$.
	Then $\Phi$ is supported on 
	\[
		\calo_F^\times \times \calo_F^\times +
			\paren{\calo_F^\times \times \calo_F^\times}j,
	\]
	with value
	\[
		\Phi\paren{\paren{x_1, x_2} + \paren{y_1, y_2}j}
			= \chi_1\paren{\frac{y_1}{x_1}} \chi_2\paren{\frac{y_2}{x_2}}.
	\]
	Again, this is invariant under $B(\calo_F)$. To find its level,
	we compute the Fourier transform of $\Phi$:
	\[
	\wh{\Phi}\paren{x_1, x_2, y_1, y_2}
		= \int_{\paren{\calo_F^\times}^4} \chi_1\paren{\frac{v_1}{u_1}}
			\chi_2\paren{\frac{v_2}{u_2}}
			\psi\paren{u_1 x_1 + u_2 x_2 - v_1 y_1 - v_2 y_2}
			du_1 du_2 dv_1 dv_2.
	\]
	This is the product of four Gaussian integrals. So we apply
	Lemma \ref{lem:epsilon}.
	$\wh{\Phi}$ is supported on 
	\[
		\Omega :=
			\varpi^{-o_1} \calo_F^\times \times \varpi^{-o_2} \calo_F^\times +
				\paren{\varpi^{-o_1}\calo_F^\times \times \varpi^{-o_2}\calo_F^\times}j,
	\]
	with value
	\[
		\wh{\Phi}\paren{\paren{x_1, x_2} + \paren{y_1, y_2}j}
			= \chi_1\paren{\frac{x_1}{y_1}} \chi_2\paren{\frac{x_2}{y_2}}
				q^{-o_1 - o_2}.
	\]
	This description shows that $\wh{\Phi}$ is invariant under
	$N(\varpi^{k}\calo_F)$ for $k = o_1 + o_2$.
	Thus we have shown that $\Phi$ is invariant under $U_0(\varpi^k)$.

	Next, we calculate
	$\wt{\Phi} = \int_{\SL_2(\calo_F)}r(k)\Phi\,dk$.
	We again use the Bruhat decomposition
	\[
		\SL_2(\calo_F)
			= w N\paren{\varpi \calo_F / \varpi^k} w U_0 \paren{\varpi^k} \, \cup \,
				N\paren{\calo_F / \varpi^k} w U_0 \paren{\varpi^k}.
	\]
	It follows that
	\[
		\wt{\Phi}
			= \paren{q^{k - 1} + q^k}^{-1}
				\paren{w \sum_{b \in \varpi \calo_F / \varpi^k}
				r\big(n(b)\big) \wh{\Phi} + 
				\sum_{b \in \calo_F / \varpi^k} r\big(n(b)\big) \wh{\Phi}}.
	\]
	As before, the two sums can be rewritten as follows.
		For $i \in \set{0, 1}$, let
		\[
			\Omega_i
				:= \set{
					\paren{u_1, u_2}
						+ \paren{v_1, v_2}j \in \Omega
					\Mid
					\frac{u_1u_2}{v_1v_2}
						\in 1 + \varpi^{k - i}\calo_F}.
		\]
		Set
		\begin{align*}
			\Psi_1
				&:= \restr{w\paren{\chi\otimes\chi^{-1}}}
					{\Omega_1}, \\
			\Psi_0
				&:= \restr{\chi\otimes\chi^{-1}}{\Omega_0}.
		\end{align*}
		Then
		\[
			\wt{\Phi}
				= q^{-k}(q + 1)^{-1}
					\paren{\Psi_1 + q\Psi_0}.
		\]
		It remains to calculate the two zeta integrals at
		$s = 1$.

		Let
		\begin{align*}
			A
				&:= \paren{
					\calo_F/\varpi^{k - 1}\calo_F}^\times, \\
			\\
			\wh{A} &:= \Hom\paren{A, \C^\times}.
		\end{align*}
		For $\omega \in \wh{A}$, inflate $\omega$ to
		$\calo_F^\times$ and extend it to $F^\times$ by
		$\omega(\varpi) = 1$. For a multiplicative character
		$\mu$ of $F^\times$, set
		\[
			I_\mu(z)
				:= \int_{\calo_F^\times}
					\mu(v)\psi(zv)dv.
		\]
		Here $dv$ is the additive Haar measure, while the
		outer integrals below use the multiplicative Haar
		measure.

		We first record an identity that treats ramified and
		unramified twists uniformly. For every multiplicative
		character $\mu$ of $F^\times$, every positive integer
		$h$, and every $u \in F^\times$, one has
		\[
			\int_{F^\times}\verts{a}
				I_\mu\paren{\varpi^{-h}a}
				I_{\mu^{-1}}
					\paren{-\varpi^{-h}au}da
				= q^{-h}\mu(u)
					\one_{\calo_F^\times}(u).
		\]
		Indeed, suppose first that $\mu$ is ramified, and set
		$r := o(\mu)$. By Lemma~\ref{lem:epsilon}, the two
		Gauss integrals can be nonzero simultaneously only
		when $\ord(a) = h - r$ and
		$u \in \calo_F^\times$. On this shell, their product is
		\[
			\mu(u)\epsilon(\mu,\psi)
				\epsilon(\mu^{-1},\psi^{-1})
				= q^{-r}\mu(u).
		\]
		Since $\verts{a} = q^{-h + r}$ and the shell has
		multiplicative volume $1$, the identity follows.

		Suppose that $\mu$ is unramified. Its restriction
		to $\calo_F^\times$ is trivial, so
		\[
			I_\mu(z)
				= J(z)
				:= \one_{\calo_F}(z)
					- q^{-1}
						\one_{\varpi^{-1}\calo_F}(z).
		\]
		After the change of variables $z = \varpi^{-h}a$,
		the left-hand side is
		\[
			q^{-h}\int_{F^\times}
				\verts{z}J(z)J(zu)dz.
		\]
		Set $j_n := J(\varpi^n)$. Thus
		$j_{-1} = -q^{-1}$,
		$j_n = 1 - q^{-1}$ for $n \geq 0$, and
		$j_n = 0$ otherwise. If $r = \ord(u)$, then a direct
		summation gives
		\[
			\sum_{n \in \Z}q^{-n}j_nj_{n + r}
				= \begin{cases}
					1 &\text{if $r = 0$,} \\
					0 &\text{if $r \neq 0$.}
				\end{cases}
		\]
		This proves the identity in the unramified case.

		For $i \in \set{1, 2}$ and
		$\omega \in \wh{A}$, put
		$\mu_{i,\omega} := \chi_i\omega$.
		Character orthogonality on $A$ gives
		\begin{align*}
			\Psi_1(x,y)
				= \frac{q^{2k}}{\verts{A}}
					\sum_{\omega \in \wh{A}}
					&I_{\mu_{1,\omega}}
						\paren{\varpi^{-o_1}x_1}
					I_{\mu_{2,\omega}}
						\paren{\varpi^{-o_2}x_2} \\
					&\cdot I_{\mu_{1,\omega}^{-1}}
						\paren{-\varpi^{-o_1}y_1}
					I_{\mu_{2,\omega}^{-1}}
						\paren{-\varpi^{-o_2}y_2}.
		\end{align*}
		Write $t_1 = (a,b) \in E^\times$ and
		$t_2 = (u,u^{-1}) \in E^1$.
		Applying the preceding identity in the $a$- and
		$b$-variables shows that the contribution of each
		$\omega \in \wh{A}$ to
		$Z\paren{\Psi_1,1}$ is
		\begin{align*}
			\frac{q^{2k}}{\verts{A}}q^{-o_1-o_2}
				\int_{\calo_F^\times}
					\mu_{1,\omega}(u)
					\mu_{2,\omega}(u^{-1})
					\xi(u,u^{-1})du
				&= \frac{q^k}{\verts{A}}.
		\end{align*}
		Indeed, $\xi = \restr{\chi^{-1}}{E^1}$ and
		$\omega(u)\omega(u^{-1}) = 1$.
		Summing over $\omega$ gives
		\[
			Z\paren{\Psi_1,1} = q^k.
		\]

		For the second term, the support of $\Psi_0$ forces
		$\ord(a) = -o_1$, $\ord(b) = -o_2$, and
		$u \in \calo_F^\times$. After multiplication by the
		factor $\xi(t_2)$ in
		Proposition~\ref{prop:rs-phi}, its character value is
		$\xi(t_2)^2$. Hence
		\[
			Z\paren{\Psi_0,1}
				= q^k\int_{\calo_E^1}\xi(t)^2dt
				= \begin{cases}
					q^k
						&\text{if $\xi^2$ is unramified}, \\
					0
						&\text{if $\xi^2$ is ramified}.
				\end{cases}
		\]
		Consequently,
		\begin{align*}
			Z\paren{\Phi,1}
				&= q^{-k}(q + 1)^{-1}
					\paren{
						Z\paren{\Psi_1,1}
						+ qZ\paren{\Psi_0,1}} \\
				&= \begin{cases}
					1
						&\text{if $\xi^2$ is unramified}, \\
					\frac{1}{q + 1}
						&\text{if $\xi^2$ is ramified}.
				\end{cases}
		\end{align*}
	\end{proof}

\subsection{Semi-ramified calculation: \texorpdfstring{$E/F$}{E/F} split}
We finish Section~\ref{sec:inner-product-new} with the last remaining
case. Assume again that $\chi = (\chi_1, \chi_2)$ is ramified,
$E/F = F \oplus F$ is split. In particular, consider the
``semi-ramified'' case wherein exactly one of
$\chi_1$ and $\chi_2$ is ramified.
Without loss of generality, we assume that $\chi_1$ is ramified and
$\chi_2$ is unramified. 

\begin{prop}
	\label{prop:new-sr}
	If $\chi = (\chi_1, \chi_2)$ with $\chi_1$
	ramified and $\chi_2$ unramified, and if
	$E = F \oplus F$, then
	\[
		\calp_\rs\paren{W^\new}
			= Z(\Phi, 1)
			= \begin{cases}
				\frac{q}{q - 1}
					&\text{if $\xi^2$ is unramified}, \\
				\frac{q}{q^2 - 1}
					&\text{if $\xi^2$ is ramified}.
			\end{cases}
	\]
\end{prop}

\begin{proof}
	Set $k := o(\chi_1)$. We take $\Phi_\chi$
	to be the restriction of $\chi_1^{-1} \otimes \one$
	to $\calo_F^\times \times \calo_F$, and set
	$\Phi := \Phi_\chi \otimes \Phi_{\chi^{-1}}$.
	Thus
	\[
		\Phi(u_1, u_2, v_1, v_2)
			= \chi_1\paren{\frac{v_1}{u_1}}
	\]
	on $(\calo_F^\times \times \calo_F)^2$, and
	$\Phi$ vanishes elsewhere. In particular, $\Phi$ is
	invariant under $B(\calo_F)$.

	We first compute its Fourier transform. Applying
	Lemma~\ref{lem:epsilon} to the first and third
	variables (cf. Tate
	\cite[Section~2.5, pp.~321--322]{tate})
	and the identity
	$\wh{\one}_{\calo_F} = \one_{\calo_F}$
	to the other two variables gives
	\[
		\wh{\Phi}(x_1, x_2, y_1, y_2)
			= q^{-k}\chi_1\paren{\frac{x_1}{y_1}}
	\]
	when $x_1, y_1 \in \varpi^{-k}\calo_F^\times$
	and $x_2, y_2 \in \calo_F$, and it vanishes
	otherwise. Hence $\wh{\Phi}$ is invariant under
	$N(\varpi^k\calo_F)$, so $\Phi$ is invariant under
	$U_0(\varpi^k)$.

	The Bruhat decomposition gives
	\[
		\SL_2(\calo_F)
			= wN\paren{\varpi\calo_F/\varpi^k\calo_F}
				wU_0(\varpi^k)
			\sqcup
			N\paren{\calo_F/\varpi^k\calo_F}
				wU_0(\varpi^k).
	\]
	For $i \in \set{0,1}$, define $\Psi_i$ by
	\[
		\Psi_i(u_1, u_2, v_1, v_2)
			:= \begin{cases}
				\chi_1\paren{u_1/v_1}
					&\text{if $u_1, v_1 \in
						\varpi^{-k}\calo_F^\times$,} \\
					&\quad\text{$u_2, v_2 \in \calo_F$, and} \\
					&\quad\text{$u_1u_2 - v_1v_2
						\in \varpi^{-i}\calo_F$,} \\
				0
					&\text{otherwise.}
			\end{cases}
	\]
	Then character orthogonality in the two Bruhat cells
	gives
	\begin{equation}
		\label{eq:wt-phi-sr}
		\wt{\Phi}
			= q^{-k}(q + 1)^{-1}
				\paren{w\Psi_1 + q\Psi_0}.
	\end{equation}

	We now calculate $w\Psi_1$. Put
	\[
		J(z)
			:= \int_{\calo_F^\times}\psi(zu)du
			= \one_{\calo_F}(z)
				-q^{-1}\one_{\varpi^{-1}\calo_F}(z).
	\]
	For $0 \leq m \leq k - 2$, set
	\[
		A_m
			:= \paren{\calo_F/
				\varpi^{k - 1 - m}\calo_F}^\times.
	\]
	Inflate every character of $A_m$ to
	$\calo_F^\times$ and extend it to $F^\times$ by
	setting $\omega(\varpi) = 1$. For
	$m + 1 \leq n \leq k - 1$, define
	\[
		f_{m, n}(x)
			:= \verts{A_m}^{-1}
				\sum_{\substack{
					\omega \in \Hom(A_m, \C^\times) \\
					o(\omega) = n - m}}
					\omega(x).
	\]
	Define the following Schwartz functions. First, let
	\[
		\Phi^{\mathrm{bd}}(x, y)
			:= \begin{cases}
				q^{2 - k}\chi_1\paren{y_1/x_1}
					&\text{if $x_1, y_1 \in \calo_F^\times$} \\
					&\quad\text{and $x_2, y_2
						\in \varpi^{1 - k}\calo_F$,} \\
				0
					&\text{otherwise.}
			\end{cases}
	\]
	For $0 \leq m \leq k - 2$, let
	\[
		\Phi_m^{\mathrm{ur}}(x, y)
			:= \begin{cases}
				q^{k - 2m}\verts{A_m}^{-1}
				\chi_1\paren{y_1/x_1}
				J(\varpi^m x_2)J(\varpi^m y_2)
					&\text{if $x_1, y_1
						\in \calo_F^\times$,} \\
				0
					&\text{otherwise.}
			\end{cases}
	\]
	Finally, let $\Phi_{m, n}$ be supported on
	$(\calo_F^\times\times
	\varpi^{-n}\calo_F^\times)^2$, with value
	\[
		\Phi_{m, n}(x, y)
			:= q^{k - m - n}\chi_1\paren{y_1/x_1}
				f_{m, n}\paren{\frac{y_1y_2}{x_1x_2}}.
	\]

	\begin{lem}
		\label{lem:newform-wtphi-split-sr}
		With the above notation,
		\[
			w\Psi_1
				= \Phi^{\mathrm{bd}}
					+\sum_{m = 0}^{k - 2}
						\paren{\Phi_m^{\mathrm{ur}}
						+\sum_{n = m + 1}^{k - 1}\Phi_{m, n}}.
		\]
	\end{lem}

	\begin{proof}[Proof of
		Lemma~\ref{lem:newform-wtphi-split-sr}]
		After replacing $u_1$ and $v_1$ by
		$\varpi^{-k}u_1$ and $\varpi^{-k}v_1$,
		respectively, the domain defining $w\Psi_1$ becomes
		\[
			D
				:= \set{(u_1, u_2, v_1, v_2)
				\in(\calo_F^\times\times\calo_F)^2
				\Mid u_1u_2 - v_1v_2
				\in \varpi^{k - 1}\calo_F}.
		\]
		Decompose $D$ as a disjoint union. The boundary
		piece is
		\[
			D_{k - 1}
				= \calo_F^\times
					\times\varpi^{k - 1}\calo_F
					\times\calo_F^\times
					\times\varpi^{k - 1}\calo_F.
		\]
		For $0 \leq m \leq k - 2$, let $D_m$ be the
		piece on which
		$u_2, v_2\in \varpi^m\calo_F^\times$.
		The condition defining $D_m$ becomes
		\[
			\frac{u_1u_2}{v_1v_2}
				\in 1 + \varpi^{k - 1 - m}\calo_F
		\]
	after dividing $u_2$ and $v_2$ by $\varpi^m$.

		On $D_{k - 1}$, the two ramified Gauss integrals
		and the two additive lattice integrals give
		$\Phi^{\mathrm{bd}}$. For $m \leq k - 2$, character
		orthogonality on $A_m$ writes the contribution
		from $D_m$ as
		\begin{align*}
			&q^{2k - 2m}\verts{A_m}^{-1}
			\sum_{\omega \in\Hom(A_m, \C^\times)}
				I_{\chi_1\omega}(\varpi^{-k}x_1)
				I_\omega(\varpi^m x_2) \\
			&\qquad\qquad\qquad\cdot
				I_{(\chi_1\omega)^{-1}}
					(-\varpi^{-k}y_1)
				I_{\omega^{-1}}(-\varpi^m y_2),
		\end{align*}
		where
		$I_\mu(z):=\int_{\calo_F^\times}
		\mu(u)\psi(zu)du$.

		The trivial character $\omega = 1$ contributes
		$\Phi_m^{\mathrm{ur}}$. Suppose that $\omega$ is ramified
		with $o(\omega) = n - m$. Since
		$o(\omega) \leq k - 1 - m < k$, the character
		$\chi_1\omega$ still has conductor exponent $k$.
		Lemma~\ref{lem:epsilon} shows that the corresponding
		product is supported on
		$(\calo_F^\times\times
		\varpi^{-n}\calo_F^\times)^2$ and equals
		\[
			q^{-k - n + m}\chi_1\paren{y_1/x_1}
				\omega\paren{\frac{y_1y_2}{x_1x_2}}.
		\]
		After multiplication by
		$q^{2k - 2m}\verts{A_m}^{-1}$ and summation over
		the characters of conductor $n - m$, this is
		$\Phi_{m, n}$. Summing the disjoint pieces proves
		the lemma.
	\end{proof}

	We apply Equation~\ref{eq:wt-phi-sr} to
	Proposition~\ref{prop:rs-phi}. On the support of each
	term in Lemma~\ref{lem:newform-wtphi-split-sr}, the
	factor $\chi_1(y_1/x_1)$ cancels the factor
	$\xi(t_2)$ in the zeta integral. First,
	\[
		Z\paren{\Phi^{\mathrm{bd}}, s}
			= \frac{q^{2 - k + (k - 1)s}}{1 - q^{-s}}.
	\]
	For $0 \leq m \leq k - 2$, one has
	\begin{align*}
		Z\paren{\Phi_m^{\mathrm{ur}}, 1}
			&= q^{k - 2m}\verts{A_m}^{-1}
				\int_{F^\times}
				\verts{b}J(\varpi^m b)^2 d^\times b \\
			&= q^{k - m}\verts{A_m}^{-1}, \\
		Z\paren{\Phi_{m, n}, 1}
			&= q^{k - m}f_{m, n}(1).
	\end{align*}
	Indeed,
	\[
		\int_{F^\times}\verts{b}
			J(\varpi^m b)^2 d^\times b
			= (1 - q^{-1})^2
				\sum_{r \geq -m} q^{-r}
				+ q^{m - 1}
			= q^m.
	\]
	The nontrivial characters of $A_m$ are partitioned
	by their conductor exponents, so
	\[
		\sum_{n = m + 1}^{k - 1}f_{m, n}(1)
			= 1 - \verts{A_m}^{-1}.
	\]
	Consequently,
	\[
		Z\paren{\Phi_m^{\mathrm{ur}}, 1}
			+\sum_{n = m + 1}^{k - 1}
				Z\paren{\Phi_{m, n}, 1}
			= q^{k - m}.
	\]

	It remains to evaluate the $\Psi_0$ term. One has
	\[
		Z\paren{q\Psi_0, s}
			= \frac{q^{1 + ks}}{1 - q^{-s}}
			\begin{cases}
				1
					&\text{if $\chi_1^2$ is unramified}, \\
				0
					&\text{if $\chi_1^2$ is ramified}.
			\end{cases}
	\]
	Combining these calculations at $s = 1$ gives
	\begin{align*}
		(q + 1)Z(\Phi,1)
			&= \frac{q^{1 - k}}{1 - q^{-1}}
				+\sum_{m = 0}^{k - 2} q^{-m} \\
			&\quad
				+\begin{cases}
					\dfrac{q}{1 - q^{-1}}
						&\text{if $\chi_1^2$
							is unramified}, \\
					0
						&\text{if $\chi_1^2$
							is ramified}.
				\end{cases}
	\end{align*}
	The first two terms sum to $(1 - q^{-1})^{-1}$.
	Since $\chi_2$ is unramified, $\xi^2$ is unramified
	if and only if $\chi_1^2$ is unramified. Therefore
	\[
		Z(\Phi,1)
			= \begin{cases}
				\dfrac{q}{q - 1}
					&\text{if $\xi^2$ is unramified}, \\
				\dfrac{q}{q^2 - 1}
					&\text{if $\xi^2$ is ramified}.
			\end{cases}
	\]
\end{proof}

\begin{remark}
	The trivial character of $A_m$ contributes
	$\Phi_m^{\mathrm{ur}}$, while the functions
	$\Phi_{m, n}$ arise from the ramified characters.
	This distinction is necessary because
	Lemma~\ref{lem:epsilon} assumes that the character
	is ramified. More explicitly, for
	$0 \leq m \leq k - 2$ and
	$m + 1 \leq n \leq k - 1$,
	\[
		f_{m, n}(1)
			= \begin{cases}
				\frac{q - 2}{q - 1}
					q^{m + 2 - k}
					&\text{if $n = m + 1$,} \\
				(q - 1)q^{n - k}
					&\text{if $n \geq m + 2$.}
			\end{cases}
	\]
	Indeed, the numbers of characters of exact conductor
	exponent $1$ and $d \geq 2$ are respectively
	$q - 2$ and $q^{d - 2}(q - 1)^2$.
	The trivial character, which would formally
	correspond to $n = m$, contributes
	$\verts{A_m}^{-1}$ through
	$\Phi_m^{\mathrm{ur}}$.
	At $s = 1$, the unramified and ramified
	contributions combine to $q^{k - m}$,
	which explains the simple final sum.
\end{remark}


\section{Proof of Theorem~\ref{thm:ratio}}
\label{sec:calculation-local}

Let $p \mid N$, and retain $K_p$ and $\chi_p$
from Section~\ref{sec:locally-dihedral}.
For the local calculation, set
$E := K_p$, $F := \Q_p$, and $\chi := \chi_p$.
Let $\xi := \chi^{1 - c_{E/F}}$.
Under the canonical isomorphism
\[
	\begin{tikzcd}[row sep = tiny]
		E^\times / F^\times
			\arrow[r, "\widesim{}"]
			& E^1 \\
		x \arrow[r, mapsto]
			& \overline{x} / x
	\end{tikzcd}
\]
the character $\xi$ corresponds to
$\restr{\chi^{-1}}{E^1}$.
When $E/F$ is a field extension,
$\varpi_E$ denotes a uniformizer of $E$.
When $E = F \oplus F$ is split,
$\varpi_1$ and $\varpi_2$ denote uniformizers
of its first and second factors, respectively.

Corollary~\ref{cor:opt-rat} gives the following:
\begin{enumerate}
	\item $\calp_\rs(W^\new) \in \Q(\xi + \xi^{-1})^\times$ and
		$\calp_\rs(W^\opt) \in \Q(\xi + \xi^{-1})$;
	\item if $\xi$ is unitary, then
		$\calp_\rs(W^\opt) \in \Q(\xi + \xi^{-1})^\times$.
\end{enumerate}

We now collect the results of Sections \ref{sec:inner-product-opt}
and \ref{sec:inner-product-new}
to determine the ratio
\[
	\sbrac{\calp_\rs\paren{W^\new} :
		\calp_\rs\paren{W^\opt}}.
\]
For the explicit calculation, we need to expand the factor from
Proposition \ref{prop:opt-ur} when $\xi$ is unramified,
\begin{align*}
	\calp_\rs\paren{W^\opt}
		&= \paren{1 + q^{-1}} L\big(\Ad(\pi), 1\big) \\
		&= \frac{q + 1}{q - 1} L\big(\xi, 1\big).
\end{align*}
In particular, if $E/F$ is inert, then
\begin{align*}
	\calp_\rs\paren{W^\opt}
		&= \frac{q^2(q + 1)}
			{(q - 1)(q - \xi(\varpi))(q - \xi(\overline{\varpi}))} \\
		&= \frac{q^2}{(q - 1)^2}.
\end{align*}
Then we obtain our main local result
by combining the above equality for unramified $\xi$
and Propositions \ref{prop:opt-ur},
\ref{prop:opt-inert}, \ref{prop:opt-split},
\ref{prop:new-ur}, \ref{prop:new-inert}, \ref{prop:new-fr},
and \ref{prop:new-sr}. This in turn gives
Theorem \ref{thm:ratio} after combining all $p$ dividing $N$.

\begin{thm}
\label{thm:comparison}
	Let $F$ be a $p$-adic field with residue field $\F_q$,
	$E/F$ be a quadratic semisimple algebra,
	$\chi$ be a character of $E^\times$, and
	$\xi$ be the ``antinorm'' $\chi^{1 - c_{E/F}}$.
	Assume the following conditions,
	\begin{enumerate}[(a)]
		\item if $p = 2$, then both $E/F$ and $\chi$ are unramified,
		\item if $E/F$ is ramified, then $\chi$ is unramified.
	\end{enumerate}
	Then the ratio
	$[\calp_\rs(W^\new) : \calp_\rs(W^\opt)]$
	is given as follows.
	\begin{enumerate}
		\item If $E/F$ is ramified, then
			\[
				\sbrac{\calp_\rs\paren{W^\new} :
					\calp_\rs\paren{W^\opt}}
					= \frac{4}{q + 1}.
			\]
		\item If $E/F$ is inert, then
			\[
				\sbrac{\calp_\rs\paren{W^\new} :
					\calp_\rs\paren{W^\opt}}
					= \begin{cases}
							1
								&\text{if $\chi$ is unramified,} \\
							\frac{(q - 1)^2}{q^2}
								&\text{if $\chi$ is ramified and $\xi$ is unramified}, \\
							\frac{\xi(-1)(q + 1)^2}{q^2}
								&\text{if $\xi$ is ramified and $\xi^2$ is unramified}, \\
							\frac{\xi(-1)(q + 1)}{q^2}
								&\text{if $\xi^2$ is ramified.}
						\end{cases}
			\]
		\item If $E = F \oplus F$ is split, then for
			$\chi = (\chi_1, \chi_2)$, 
			\[
				\sbrac{\calp_\rs\paren{W^\new} :
					\calp_\rs\paren{W^\opt}}
					= \begin{cases}
							1
								&\text{if $\chi$ is unramified,} \\
							\frac{(q - 1)(q - \xi(\varpi_1))
							(q - \xi(\varpi_2))}{(q + 1) q^2}
								&\text{if $\chi$ is ramified and $\xi$ is unramified,} \\
							\frac{\chi(-1) (q - 1)^2 q^{2o(\xi) - 1}}
							{q^{2o(\xi) + 1} - q^{2o(\xi)} + 2}
								&\text{if $\xi$ is ramified, $\xi^2$ is unramified,} \\
								&\text{and exactly one of the $\chi_i$ is ramified,} \\
							\frac{\chi(-1) (q - 1)^3 q^{2o(\xi) - 2}}
							{q^{2o(\xi) + 1} - q^{2o(\xi)} + 2}
								&\text{if $\xi$ is ramified, $\xi^2$ is unramified,} \\
								&\text{and both $\chi_i$ are ramified,} \\
							\frac{\chi(-1) (q - 1)^2 q^{2o(\xi) - 1}}
							{(q + 1) \paren{q^{2o(\xi) + 1} - q^{2o(\xi)} + 2}}
								&\text{if $\xi^2$ is ramified} \\
								&\text{and exactly one of the $\chi_i$ is ramified,} \\
							\frac{\chi(-1) (q - 1)^3 q^{2o(\xi) - 2}}
							{(q + 1)(q^{2o(\xi) + 1} - q^{2o(\xi)} + 2)}
								&\text{if $\xi^2$, $\chi_1$, and $\chi_2$ are all ramified.}
						\end{cases}
			\]
	\end{enumerate}
\end{thm}


\section{The constants
	\texorpdfstring{$c_{f, \rs}$}{c f,RS} and
	\texorpdfstring{$m_f$}{m f}}
\label{sec:petersson}

\subsection{Petersson norms of newforms}

Let $f$ be a newform of weight $1$ for $\Gamma_0(N)$ for some $N$
with central character $\omega$. 
Let $f^*$ be the dual form
$f^*(z) := \overline{f(-\overline{z})}$.
Equivalently, if $f$ has the $q$-expansion
$f = \sum_n a_nq^n$,
then $f^*$ has the $q$-expansion 
$f^*= \sum_n \overline{a_n} q^n$ with complex conjugate coefficients.
Let $d\mu$ be the $\SL_2(\Z)$-invariant hyperbolic volume form
$\frac{dx dy}{y^2}$ on the upper half-plane.
The main purpose of this section is to evaluate
the Petersson norm,
\[
	\norm{f}_\R^2 = \int_{X_0(N)} \verts{f(z)}^2 y \, d\mu.
\]

Stark \cite[Section~6]{stark-1975}
observed that
$\norm{f}_\R^2$ is a multiple of
$L(\Ad(\rho_f), 1)$
by finitely many Euler factors,
and gave an explicit calculation for the
weight-$1$ theta series attached to
$\Q(\sqrt{-23})$.
We first record the relation between this $L$-value
and $L'(\Ad(\rho_f), 0)$.
Let $N_{\Ad}$ be the Artin conductor of
$\Ad(\rho_f)$ and set
$\Gamma_\R(s) := \pi^{-s/2}\Gamma(s/2)$.
Since complex conjugation acts on $\Ad(\rho_f)$
with eigenvalues $1$, $-1$, and $-1$, its completed
Artin $L$-function is
\[
	\Lambda_{\Ad}(s)
		:= N_{\Ad}^{s/2}\Gamma_\R(s)
			\Gamma_\R(s + 1)^2
			L\paren{\Ad(\rho_f), s}.
\]
Since $\Ad(\rho_f)$ is self-dual, its functional
equation \cite[Section~1]{stark-1975} is
$\Lambda_{\Ad}(s)
= \epsilon_{\Ad}\Lambda_{\Ad}(1 - s)$,
where $\epsilon_{\Ad}$ is its root number.
Using the simple zero of $L(\Ad(\rho_f), s)$
at $s = 0$, together with
$\Gamma_\R(s) \sim 2/s$,
$\Gamma_\R(1) = 1$, and
$\Gamma_\R(2) = \pi^{-1}$, gives
\[
	L\paren{\Ad(\rho_f), 1}
		= 2\epsilon_{\Ad}^{-1}\pi^2
			N_{\Ad}^{-1/2}
			L'\paren{\Ad(\rho_f), 0}.
\]

The point of this section is to give a uniform formula
for the Rankin--Selberg constant in the normalization
used here and to identify its ramified local factors
with the comparison constants calculated above.
Proposition~\ref{prop:crs} expresses the algebraic
constant $c_{f, \rs}$ in
\[
	\norm{f}_\R^2
		= c_{f, \rs}
			L'\paren{\Ad(\rho_f), 0}
\]
as a product of ramified local Euler-factor ratios.
Together with the preceding functional equation, this
also determines the proportionality constant in
\[
	\norm{f}_\R^2
		= \frac{\epsilon_{\Ad}N_{\Ad}^{1/2}}
			{2\pi^2}c_{f,\rs}
			L\paren{\Ad(\rho_f),1}.
\]
For imaginary dihedral forms,
Corollary~\ref{cor:crs} combines this formula with
Theorem~\ref{thm:ratio} to evaluate $c_{f,\rs}$ in
terms of the local comparison constants
$\alpha_{\chi_p}$.

Let
$\zeta^*(s) := \pi^{-s/2}\Gamma(s/2)\zeta(s)$
be the completed Riemann zeta function. Our approach uses
the completed Eisenstein series
\[
	E(z, s)
		= \frac{\pi^{-s}\Gamma(s)}{2}
			\sum_{(m, n) \in\Z^2\setminus\set{(0, 0)}}
			\frac{y^s}{\verts{mz + n}^{2s}}
		= \frac{\zeta^*(2s)}{2}
			\sum_{\gcd(m, n) = 1}
			\frac{y^s}{\verts{mz + n}^{2s}}.
\]
Jacquet--Zagier \cite[Section~1]{jacquet-zagier}
show that $E(z, s)$ has simple poles at $s = 1$ and
$s = 0$ with respective residues $1/2$ and $-1/2$.
Define the Rankin--Selberg convolution by
\[
	\Lambda\paren{f\otimes f^*, s}
		:= \int_{X_0(N)}\verts{f(z)}^2y E(z, s) d\mu.
\]
Since $f$ is cuspidal, the integral is rapidly decreasing
at the cusps, so we may take the residue under the
integral. The residue of $E(z, s)$ at $s = 0$ therefore gives
\[
	\Res_{s = 0}\Lambda\paren{f\otimes f^*, s}
		= -\frac{1}{2}\norm{f}_\R^2,
\]
or equivalently
\[
	\norm{f}_\R^2
		= -2\Res_{s = 0}
			\Lambda\paren{f\otimes f^*, s}.
\]

To evaluate the Rankin--Selberg convolution, define
\[
	\Tr_{\Gamma_0(N)}^{\SL_2(\Z)}\paren{\verts{f}^2}
		:= \sum_{\gamma \in\Gamma_0(N)\bs\SL_2(\Z)}
			\verts{\restr{f}{1}\gamma}^2.
\]
Then unfolding gives
\begin{align*}
	\Lambda\paren{f\otimes f^*, s}
		&= \int_{X(1)}
			\Tr_{\Gamma_0(N)}^{\SL_2(\Z)}
				\paren{\verts{f}^2}(z)
			y E(z, s) d\mu \\
		&= \zeta^*(2s)
			\int_{\Gamma_\infty\bs\calh}
			\Tr_{\Gamma_0(N)}^{\SL_2(\Z)}
				\paren{\verts{f}^2}(z)
			y^{1 + s} d\mu \\
		&= \zeta^*(2s)\int_0^\infty
			\paren{\int_0^1
				\Tr_{\Gamma_0(N)}^{\SL_2(\Z)}
					\paren{\verts{f}^2}(x + yi) dx}
			y^s\frac{dy}{y}.
\end{align*}
For the inner integral, use the $q$-expansion
$\restr{f}{\gamma}(x + iy)
	= \sum_n a_{\gamma,n}
		e^{2\pi i n(x + iy)}$.
It follows that
\[
	\int_0^1
		\Tr\paren{\verts{f(x + iy)}^2}dx
			= \sum_n b_n e^{-4\pi n y},
\]
where
\[
	b_n
		:= \sum_{\gamma \in
			\Gamma_0(N) \bs \SL_2(\Z)}
			\verts{a_{\gamma,n}}^2.
\]
Therefore
\begin{align*}
	\Lambda\paren{f \otimes f^*, s}
		&= \zeta^*(2s)
			\int_0^\infty
				\sum_n b_n e^{-4\pi n y}
				y^s\frac{dy}{y} \\
		&= \zeta^*(2s)
			\sum_n b_n
				\int_0^\infty
					e^{-4\pi n y}
					y^s\frac{dy}{y} \\
		&= (2\pi)^{-2s}\Gamma(s)^2
			\zeta(2s)\sum_n\frac{b_n}{n^s}.
\end{align*}
where in the last equality, we use the substitution
$y \mapsto y/(4 \pi n)$.

Write $L(f \times f^*, s) = \zeta(2s) \sum_n \frac{b_n}{n^s}$.
As $b_n$ is multiplicative with a prime decomposition,
so is $L(f \times f^*, s)$:
\[
	L\paren{f \times f^*, s}
		= \prod_p L_p\paren{f \times f^*, s},
\]
where the factor for each prime $p$ is
\[
	L_p\paren{f \times f^*, s}
		= \zeta_p(2s) \sum_{n = 0}^\infty \frac{b_{p^n}}{p^{ns}}.
\]
We now demonstrate Proposition \ref{prop:crs},
the proportionality of the Petersson norm $\norm{f}_{\R}^2$
to $L'(\Ad(\rho_f), 0)$ by
\[
	c_{f, \rs} = \restr{\paren{\prod_{p \mid N}
		\frac{L_p\paren{f \times f^*, s}}
		{\zeta_p(s) L_p\big(\Ad(\rho_f), s\big)}}}
		{s = 0} \in \Q(\chi_{\rho_f}).
\]

\begin{proof}[Proof of Proposition \ref{prop:crs}]
	If $p \nmid N$, then
	\[
		L_p(f, s)
			= \sum_{n = 0}^\infty \frac{a_{p^n}}{p^{ns}}
			= \frac{1}{1 - a_p p^{-s} + \omega(p) p^{-2s}}.
	\]
	Take the decomposition,
	\[
		1 - a_p p^{-s} + \omega(p)p^{-2s}
			= \paren{1 - \alpha_{p, 1} p^{-s}}\paren{1 - \alpha_{p, 2} p^{-s}}.
	\]
	Then
	\[
		\frac{1}{1 - a_p p^{-s} + \omega(p)p^{-2s}}
			= \paren{\sum_\ell \alpha_{p, 1}^\ell p^{-\ell s}}
				\paren{\sum_m \alpha_{p, 2}^m p^{-ms}}.
	\]
	It follows that
	\[
		a_{p^n}
			= \sum_{\ell + m = n} \alpha_{p, 1}^\ell \alpha_{p, 2}^m
			= \frac{\alpha_{p, 1}^{n + 1} - \alpha_{p, 2}^{n + 1}}{\alpha_{p, 1} - \alpha_{p, 2}}.
	\]
	Notice that $\alpha_{p, 1}$ and $\alpha_{p, 2}$ are
	eigenvalues of representation of finite groups.
	Thus $\alpha_{p, 1}$ and $\alpha_{p, 2}$ have norm $1$.
	This shows that $f^*$ has eigenvalues
	$\beta_{p, 1} := \alpha_{p, 1}^{-1}$ and
	$\beta_{p, 2} := \alpha_{p, 2}^{-1}$.
	By the calculation in Example \ref{ex:rs-unramified},
	\[
		\sum_{n=0}^\infty \frac{b_{p^n}}{p^{ns}}
			= \sum_{n=0}^\infty
				\frac{\alpha_{p, 1}^{n + 1} - \alpha_{p, 2}^{n + 1}}{\alpha_{p, 1} - \alpha_{p, 2}}
				\frac{\beta_{p, 1}^{n + 1} - \beta_{p, 2}^{n + 1}}{\beta_{p, 1} - \beta_{p, 2}}
				p^{-ns}
			= \paren{1 + p^{-s}} L_p\big(\Ad(\rho_f), s\big).
	\]
	Therefore for $p \nmid N$,
	\[
		L_p\paren{f \times f^*, s}
			= \zeta_p(2s) \paren{1 + p^{-s}} L_p\big(\Ad(\rho_f), s\big)
			= \zeta_p(s) L_p\big(\Ad(\rho_f), s\big).
	\]

	In summary, we have shown that
	\begin{align*}
		L\paren{f \times f^*,s}
			&= C(s)\zeta(s)
				L\paren{\Ad(\rho_f),s}, \\
		\Lambda\paren{f \times f^*,s}
			&= C(s)(2\pi)^{-2s}\Gamma(s)^2
				\zeta(s)L\paren{\Ad(\rho_f),s},
	\end{align*}
	where $C(s)$ is a product of polynomials in
	$p^{-s}$ for $p \mid N$, with coefficients in
	$\Q(\chi_{\rho_f})$:
	\[
		C(s)
			:= \prod_{p \mid N}
				\frac{
					L_p\paren{f \times f^*,s}}
					{\zeta_p(s)
					L_p\paren{\Ad(\rho_f),s}}.
	\]
	Notice that $\zeta(0) = -1/2$, and that $\Gamma(s)$ has a simple pole
	at $s = 0$ with residue $1$. Then take
	$\norm{f}_\R^2 = -2 \, \Res_{s = 0}
	\Lambda\paren{f \otimes f^*, s}$ to obtain the desired relation.
\end{proof}

\subsection{Proof of Corollary~\ref{cor:crs}}
Let $f$ be imaginary dihedral and set
$f^\new := f \otimes f^*$.
By Theorem~\ref{thm:zhang-hvs} and the normalization
of the Rankin--Selberg period,
\[
	\calp_\rs\paren{f^\new}
		= \sbrac{\PSL_2(\Z) : \Gamma_0(N)}^{-1}
			\Reg_\R\paren{u_f}.
\]
Set
\[
	u_{\varphi^\opt}
		:= -\frac{\sbrac{H_c : H_1}w_K}
			{2\sbrac{\PSL_2(\Z) : \Gamma_0(N)}}
			u_\stark.
\]
By \cite[Lemma~5.2]{zhang-hvs},
\[
	\calp_\rs\paren{\varphi^\opt}
		= \Reg_\R\paren{u_{\varphi^\opt}}
		= -\frac{\sbrac{H_c : H_1}w_K}
			{2\sbrac{\PSL_2(\Z) : \Gamma_0(N)}}
			\Reg_\R\paren{u_\stark}.
\]
Therefore
\[
	c_{f, \rs}
		= -\frac{\sbrac{H_c : H_1}w_K}{2}
			\frac{\calp_\rs\paren{f^\new}}
			{\calp_\rs\paren{\varphi^\opt}}.
\]
Theorem~\ref{thm:ratio} proves the formula in
Corollary~\ref{cor:crs}. The uniqueness in
Theorem~\ref{thm:zhang-hvs} also gives
\begin{equation}
\label{eq:canonical-unit}
	u_f
		= c_{f, \rs}u_\stark
		= -\frac{\sbrac{H_c : H_1}w_K}{2}
			\paren{\prod_{q \mid N}\alpha_{\chi_q}}
			u_\stark
\end{equation}
in $\calu\paren{\Ad(\rho_f)} \otimes_\Z \Q$.


\subsection{Proof of Corollary~\ref{cor:m}}
\label{sec:cor-m}
For a prime $p \nmid 6N$, let
$\Delta_p := \F_p^\times \otimes_\Z \Z\sbrac{1/6}$.
The prime-to-$6$ quotient of the Shimura covering
$X_1(p) \lra X_0(p)$ defines a class
\[
	\mfs_p
		\in H^1_{\mathrm{\acute{e}t}}
			\paren{
				X_0(p)_{\Z[1/(6p)]},
				\Delta_p}.
\]
By Serre duality, it gives a functional on weight-$2$
forms, also denoted by $\mfs_p$; see
\cite[Section~3.1]{hv} and
\cite[Sections~2.2--2.4]{zhang-hv}.
For $\ell^t \Vert p - 1$ with $\ell \geq 5$ and
$\ell \nmid N$, fix a discrete logarithm
$\log_\ell: \Delta_p \lra \Z/\ell^t\Z$, and set
$\calp_{\hv, \ell}
	:= \paren{1 \otimes \log_\ell}
		\circ\calp_\hv$.
The finite-level comparison is
\begin{equation}
	\label{eq:shimura-newform}
	\log_\ell\norm{f}_{\F_p^\times}^2
		= \sbrac{\PSL_2(\Z) : \Gamma_0(N)}
			\calp_{\hv, \ell}\paren{f^\new}.
\end{equation}

We now apply the integral comparison of
\cite[Section~4]{zhang-hv},
in which we have
the imaginary quadratic field $K$ and the
finite character $\chi$ associated to $f$,
so $\rho_f \cong \Ind_{G_K}^{G_\Q}\paren{\chi}$.
Let $\xi := \chi^{1 - c_K}$, and let
$L := \Q(\chi)$ be the field generated by the values
of $\chi$.

The character of $\rho_f$ is zero outside $G_K$ and
is $\chi+\chi^{c_K}$ on $G_K$. Hence
$\Q(f) \subseteq L$.
Set
$F_f := \Q\paren{\chi_{\Ad(\rho_f)}}$.
Since
$\chi_{\Ad(\rho_f)}
= \chi_{\rho_f}\chi_{\rho_f^\vee}-1$,
one also has $F_f \subseteq \Q(f)$.
Thus $F_f \subseteq \Q(f) \subseteq L$.
If $\zeta$ generates $\Im(\xi)$, then
$F_f = \Q(\zeta + \zeta^{-1})$ and
$R_f = \Z[\zeta + \zeta^{-1}]
= \calo_{F_f}$.
In particular,
$\Tr_{L/F_f}(\calo_L) \subseteq R_f$, where
$\calo_L$ is the ring of integers of $L$.

The newform integral model is defined over
$\Q(f)$, while the optimal vector
$\varphi^\opt$ and its associated unit
$u_{\varphi^\opt}$ are naturally defined over $L$;
see \cite[Section~9]{zhang-hv}.
We therefore base change the integral model to
$\calo_L\sbrac{1/N}$ and define
\[
	\calu_L\paren{\Ad(\rho_f)}
		:= \calu\paren{\Ad(\rho_f)}
			\otimes_{R_f}\calo_L.
\]
Let $\Sigma$ be the set of primes dividing $N$, and
let
$\Pi_{f,\Sigma,\calo_L[1/N]}$
be the integral model used in
\cite[Section~4]{zhang-hv}.
Choose the least positive integer $c_L^\opt$ such that
\[
	\varphi_L^\opt
		:= c_L^\opt\varphi^\opt
		\in \Pi_{f,\Sigma,\calo_L[1/N]},
\]
and let $f_{\varphi_L^\opt}$ be its modular avatar.
Choose the least positive integer $d_L^\opt$ such that
$d_L^\opt u_{\varphi^\opt}$ belongs to
$\calu_L\paren{\Ad(\rho_f)}$.
By \cite[Lemma~5.2]{zhang-hvs},
\begin{equation}
	\label{eq:integral-optimal-period}
	d_L^\opt
		\calp_{\hv, \ell}
			\paren{f_{\varphi_L^\opt}}
		= \log_\ell\Reg_{\F_p^\times}
			\paren{
				c_L^\opt d_L^\opt
				u_{\varphi^\opt}}.
\end{equation}

Since the Rankin--Selberg period is a nonzero
multiple of the canonical pairing $\calp$,
Theorem~\ref{thm:ratio} gives
\[
	\frac{\calp\paren{\varphi_L^\opt}}
		{\calp\paren{f^\new}}
		= \frac{c_L^\opt}
			{\prod_{q \mid N}\alpha_{\chi_q}}.
\]
Choose nonzero $a_f, b_f \in \calo_L$ such that this
ratio is $a_f/b_f$ and such that
$\verts{\rmn_{L/\Q}\paren{a_f}}$ is minimal among
these choices. Set
\[
	a_f^\vee
		:= \frac{\verts{\rmn_{L/\Q}\paren{a_f}}}{a_f}
		\in \calo_L.
\]
Thus $a_f^\vee a_f = \verts{\rmn_{L/\Q}\paren{a_f}}$.

Define
\begin{equation}
	\label{eq:explicit-multiplier}
	m_f
		:= \sbrac{L : F_f}
			\verts{\rmn_{L/\Q}\paren{a_f}}
			d_L^\opt
			\paren{\prod_{q \mid N}(q - 1)}^2.
\end{equation}
The integral comparison in
\cite[Section~4]{zhang-hv} gives
\[
	b_f\paren{\prod_{q \mid N}(q - 1)}^2
		\calp_{\hv, \ell}
			\paren{f_{\varphi_L^\opt}}
		= a_f\paren{\prod_{q \mid N}(q - 1)}^2
			\calp_{\hv, \ell}\paren{f^\new}.
\]
Combining this with
Equation~\ref{eq:integral-optimal-period}, multiplying
by $a_f^\vee$, taking the trace from $L$ to $F_f$,
and using Equation~\ref{eq:shimura-newform} gives
\begin{align*}
	m_f\log_\ell\norm{f}_{\F_p^\times}^2
		&= \log_\ell\Reg_{\F_p^\times}
		\Bigg(
			\sbrac{\PSL_2(\Z) : \Gamma_0(N)}
			\Tr_{L/F_f} \paren{
				a_f^\vee b_f
				\paren{\prod_{q \mid N}(q - 1)}^2
				c_L^\opt d_L^\opt
				u_{\varphi^\opt}}
		\Bigg).
\end{align*}

The definitions of $a_f/b_f$ and $u_f$ give
\begin{align*}
	b_f c_L^\opt
		&= a_f\prod_{q \mid N}\alpha_{\chi_q}, \\
	u_f
		&= \sbrac{\PSL_2(\Z) : \Gamma_0(N)}
			\paren{\prod_{q \mid N}\alpha_{\chi_q}}
			u_{\varphi^\opt}.
\end{align*}
Consequently, the unit in the preceding regulator is
\[
	\Tr_{L/F_f}
		\paren{
			a_f^\vee a_f
			\paren{\prod_{q \mid N}(q - 1)}^2
			d_L^\opt u_f}.
\]
Since
$a_f^\vee a_f
= \verts{\rmn_{L/\Q}\paren{a_f}}$
and $u_f$ is defined over $F_f$, this unit is
\[
	\sbrac{L : F_f}
		\verts{\rmn_{L/\Q}\paren{a_f}}
		\paren{\prod_{q \mid N}(q - 1)}^2
		d_L^\opt u_f
		= m_f u_f.
\]
Finally,
$\Tr_{L/F_f}(\calo_L)\subseteq R_f$, so
$m_f u_f$ belongs to
$\calu\paren{\Ad(\rho_f)}$.
We have therefore proved
\[
	m_f\log_\ell\norm{f}_{\F_p^\times}^2
		= \log_\ell\Reg_{\F_p^\times}
			\paren{m_f u_f}
\]
for every $p$ and $\ell$ fixed above.
Reassembling the primary components gives
\[
	m_f\norm{f}_{\F_p^\times}^2
		= \Reg_{\F_p^\times}\paren{m_f u_f}
\]
for every prime $p \nmid 6N$.
Theorem~\ref{thm:zhang-hvs} and the linearity of the
Stark regulator also give
$\Reg_\R\paren{m_f u_f} = m_f\norm{f}_\R^2$.


\bibliography{bibliography}{}

\begin{thebibliography}{DHRV22}

\bibitem[Bum97]{bump}
Daniel Bump.
\newblock {\em Automorphic forms and representations}, volume~55 of {\em
  Cambridge Studies in Advanced Mathematics}.
\newblock Cambridge University Press, Cambridge, 1997.

\bibitem[Cas73]{casselman}
William Casselman.
\newblock On some results of {A}tkin and {L}ehner.
\newblock {\em Math. Ann.}, 201:301--314, 1973.

\bibitem[DHRV22]{dhrv}
Henri Darmon, Michael Harris, Victor Rotger, and Akshay Venkatesh.
\newblock The derived {H}ecke algebra for dihedral weight one forms.
\newblock {\em Michigan Math. J.}, 72:145--207, 2022.

\bibitem[HK91]{harris-kudla-1991}
Michael Harris and Stephen~S. Kudla.
\newblock The central critical value of a triple product {$L$}-function.
\newblock {\em Ann. of Math. (2)}, 133(3):605--672, 1991.

\bibitem[HK04]{harris-kudla-2004}
Michael Harris and Stephen~S. Kudla.
\newblock On a conjecture of {J}acquet.
\newblock In {\em Contributions to automorphic forms, geometry, and number
  theory}, pages 355--371. Johns Hopkins Univ. Press, Baltimore, MD, 2004.

\bibitem[Hor23]{horawa}
Aleksander Horawa.
\newblock Motivic action on coherent cohomology of {H}ilbert modular varieties.
\newblock {\em Int. Math. Res. Not. IMRN}, (12):10439--10531, 2023.

\bibitem[HV19]{hv}
Michael Harris and Akshay Venkatesh.
\newblock Derived {H}ecke algebra for weight one forms.
\newblock {\em Exp. Math.}, 28(3):342--361, 2019.

\bibitem[Jac72]{jacquet}
Herv\'{e} Jacquet.
\newblock {\em Automorphic forms on {${\rm GL}(2)$}. {P}art {II}}.
\newblock Lecture Notes in Mathematics, Vol. 278. Springer-Verlag, Berlin-New
  York, 1972.

\bibitem[JL70]{jacquet-langlands}
Herv\'{e} Jacquet and Robert~P. Langlands.
\newblock {\em Automorphic forms on {${\rm GL}(2)$}}.
\newblock Lecture Notes in Mathematics, Vol. 114. Springer-Verlag, Berlin-New
  York, 1970.

\bibitem[JZ87]{jacquet-zagier}
Herv\'{e} Jacquet and Don Zagier.
\newblock Eisenstein series and the {S}elberg trace formula. {II}.
\newblock {\em Trans. Amer. Math. Soc.}, 300(1):1--48, 1987.

\bibitem[Lec23]{lecouturier-hv}
Emmanuel Lecouturier.
\newblock On triple product {L}-functions and a conjecture of
  {H}arris-{V}enkatesh.
\newblock {\em Int. Math. Res. Not. IMRN}, (22):19476--19506, 2023.

\bibitem[Mar21]{marcil}
David Marcil.
\newblock Numerical verification of a conjecture of {H}arris and {V}enkatesh.
\newblock {\em J. Number Theory}, 221:484--495, 2021.

\bibitem[Sta71]{stark-1971}
Harold~M. Stark.
\newblock Values of {$L$}-functions at {$s=1$}. {I}. {$L$}-functions for
  quadratic forms.
\newblock {\em Adv. Math.}, 7:301--343 (1971), 1971.

\bibitem[Sta75]{stark-1975}
Harold~M. Stark.
\newblock {$L$}-functions at {$s=1$}. {II}. {A}rtin {$L$}-functions with
  rational characters.
\newblock {\em Adv. Math.}, 17(1):60--92, 1975.

\bibitem[Sta76]{stark-1976}
Harold~M. Stark.
\newblock {$L$}-functions at {$s=1$}. {III}. {T}otally real fields and
  {H}ilbert's twelfth problem.
\newblock {\em Adv. Math.}, 22(1):64--84, 1976.

\bibitem[Sta80]{stark-1980}
Harold~M. Stark.
\newblock {$L$}-functions at {$s=1$}. {IV}. {F}irst derivatives at {$s=0$}.
\newblock {\em Adv. Math.}, 35(3):197--235, 1980.

\bibitem[Tat67]{tate}
John~T. Tate.
\newblock Fourier analysis in number fields, and {H}ecke's zeta-functions.
\newblock In {\em Algebraic {N}umber {T}heory ({P}roc. {I}nstructional {C}onf.,
  {B}righton, 1965)}, pages 305--347. Academic Press, London, 1967.

\bibitem[Wal85]{waldspurger}
Jean-Loup Waldspurger.
\newblock Sur les valeurs de certaines fonctions {$L$} automorphes en leur
  centre de sym\'{e}trie.
\newblock {\em Compositio Math.}, 54(2):173--242, 1985.

\bibitem[Wei64]{weil}
Andr\'{e} Weil.
\newblock Sur certains groupes d'op\'{e}rateurs unitaires.
\newblock {\em Acta Math.}, 111:143--211, 1964.

\bibitem[YZZ13]{yzz}
Xinyi Yuan, Shou-Wu Zhang, and Wei Zhang.
\newblock {\em The {G}ross-{Z}agier formula on {S}himura curves}, volume 184 of
  {\em Annals of Mathematics Studies}.
\newblock Princeton University Press, Princeton, NJ, 2013.

\bibitem[Zha01]{zhang-asian}
Shou-Wu Zhang.
\newblock Gross-{Z}agier formula for {${\rm GL}_2$}.
\newblock {\em Asian J. Math.}, 5(2):183--290, 2001.

\bibitem[Zha23a]{zhang-hv}
Robin Zhang.
\newblock The {H}arris--{V}enkatesh conjecture for derived {H}ecke operators
  {I}: imaginary dihedral forms, 2023.
\newblock arXiv:2301.00570.

\bibitem[Zha23b]{zhang-hvs}
Robin Zhang.
\newblock The {H}arris--{V}enkatesh conjecture for derived {H}ecke operators
  {II}: a unified {S}tark conjecture, 2023.
\newblock arXiv:2305.08956.

\bibitem[Zha23c]{zhang-real}
Robin Zhang.
\newblock The {H}arris--{V}enkatesh conjecture for derived {H}ecke operators
  {IV}: real dihedral forms, 2023.
\newblock In progress.

\end{thebibliography}
\bibliographystyle{alpha}

\end{document}